\documentclass[11pt]{article}
\usepackage{latexsym,amssymb,amsmath}
\usepackage{amssymb,amsthm,upref,amscd}
\usepackage{color}
\usepackage[T1]{fontenc}
\begin{document}
\baselineskip=15pt

\numberwithin{equation}{section}

\pretolerance1000
\newtheorem{theorem}{Theorem}[section]
\newtheorem{lemma}[theorem]{Lemma}
\newtheorem{proposition}[theorem]{Proposition}
\newtheorem{corollary}[theorem]{Corollary}
\newtheorem{remark}[theorem]{Remark}
\newtheorem{definition}[theorem]{Definition}
\renewcommand{\theequation}{\thesection.\arabic{equation}}

\newcommand{\A}{\mathbb{A}}
\newcommand{\B}{\mathbb{B}}
\newcommand{\C}{\mathbb{C}}
\newcommand{\D}{\mathbb{D}}
\newcommand{\E}{\mathbb{E}}
\newcommand{\F}{\mathbb{F}}
\newcommand{\G}{\mathbb{G}}
\newcommand{\I}{\mathbb{I}}
\newcommand{\J}{\mathbb{J}}
\newcommand{\K}{\mathbb{K}}
\newcommand{\M}{\mathbb{M}}
\newcommand{\N}{\mathbb{N}}
\newcommand{\Q}{\mathbb{Q}}
\newcommand{\R}{\mathbb{R}}
\newcommand{\T}{\mathbb{T}}
\newcommand{\U}{\mathbb{U}}
\newcommand{\V}{\mathbb{V}}
\newcommand{\W}{\mathbb{W}}
\newcommand{\X}{\mathbb{X}}
\newcommand{\Y}{\mathbb{Y}}
\newcommand{\Z}{\mathbb{Z}}
\newcommand\ca{\mathcal{A}}
\newcommand\cb{\mathcal{B}}
\newcommand\cc{\mathcal{C}}
\newcommand\cd{\mathcal{D}}
\newcommand\ce{\mathcal{E}}
\newcommand\cf{\mathcal{F}}
\newcommand\cg{\mathcal{G}}
\newcommand\ch{\mathcal{H}}
\newcommand\ci{\mathcal{I}}
\newcommand\cj{\mathcal{J}}
\newcommand\ck{\mathcal{K}}
\newcommand\cl{\mathcal{L}}
\newcommand\cm{\mathcal{M}}
\newcommand\cn{\mathcal{N}}
\newcommand\co{\mathcal{O}}
\newcommand\cp{\mathcal{P}}
\newcommand\cq{\mathcal{Q}}
\newcommand\rr{\mathcal{R}}
\newcommand\cs{\mathcal{S}}
\newcommand\ct{\mathcal{T}}
\newcommand\cu{\mathcal{U}}
\newcommand\cv{\mathcal{V}}
\newcommand\cw{\mathcal{W}}
\newcommand\cx{\mathcal{X}}
\newcommand\ocd{\overline{\cd}}

\def\c{\centerline}
\def\ov{\overline}
\def\emp {\emptyset}
\def\pa {\partial}
\def\bl{\setminus}
\def\op{\oplus}
\def\sbt{\subset}
\def\un{\underline}
\def\al {\alpha}
\def\bt {\beta}
\def\de {\delta}
\def\Ga {\Gamma}
\def\ga {\gamma}
\def\lm {\lambda}
\def\Lam {\Lambda}
\def\om {\omega}
\def\Om {\Omega}
\def\sa {\sigma}
\def\vr {\varepsilon}
\def\va {\varphi}

\title{\bf Existence of positive multi-bump solutions for a  Schr\"odinger-Poisson system in $\mathbb{R}^{3}$ }

\author{Claudianor O. Alves$^a$\thanks{C. O. Alves was partially supported by CNPq/Brazil
 301807/2013-2 and INCT-MAT, coalves@dme.ufcg.edu.br}\, , \, \   Minbo Yang$^b$\thanks{M. Yang was supported by NSFC (11101374, 11271331), mbyang@zjnu.edu.cn}\vspace{2mm}
\and {\small $a.$ Universidade Federal de Campina Grande} \\ {\small Unidade Acad\^emica de Matem\'{a}tica} \\ {\small CEP: 58429-900, Campina Grande - Pb, Brazil}\\
{\small $b.$ Department of Mathematics, Zhejiang Normal University} \\ {\small Jinhua, 321004, P. R. China.}}

\date{}
\maketitle

\begin{abstract}
In this paper we are going to study a class of Schr\"odinger-Poisson system
$$
\left\{ \begin{array}{ll}
  - \Delta u + (\lambda a(x)+1)u+ \phi u = f(u)  \mbox{ in } \,\,\, \mathbb{R}^{3},\\
-\Delta \phi=u^2  \mbox{ in } \,\,\, \mathbb{R}^{3}.\\
\end{array}
\right.
$$
Assuming that the nonnegative function $a(x)$ has a potential well $int (a^{-1}(\{0\}))$  consisting of $k$ disjoint components $\Omega_1, \Omega_2, ....., \Omega_k$ and the nonlinearity $f(t)$ has a subcritical growth, we are able to establish the existence of positive multi-bump solutions by variational methods.

 \vspace{0.3cm}


\vspace{0.3cm}

\end{abstract}

\section{Introduction}

This paper was motivated by some recent works concerning the nonlinear Schr\"{o}dinger-Poisson
system
$$
\left\{ \begin{array}{l}
 -i\frac{\partial\psi}{\partial t} = - \Delta \psi + V(x)\psi +
\phi(x)\psi - |\psi|^{p-2}\psi  \mbox{ in $\mathbb{R}^{3}$},\\
-\Delta \phi=|\psi|^2  \mbox{ in $\mathbb{R}^{3}$},\\
\end{array}
\right.\eqno{(NSP)}
$$
where $V:\mathbb{R}^{3} \to \mathbb{R}$ is a nonnegative continuous function with $\displaystyle \inf_{x \in \mathbb{R}^{3}}V(x)>0$, $2<p< 2^*=6$ and $\psi : \overline{\Omega} \to \mathbb{C}$ and $\phi: \overline{\Omega} \to \mathbb{R}$ are two unknown functions.

The first equation in system $(NSP)$, called Schr\"odinger equation,
describes quantum (non-relativistic) particles interacting with the
eletromagnetic field generated by the motion.   An interesting phenomenon of
this Schr\"odinger type equation is that the potential $\phi(x)$ is
determined by the charge of wave function itself, that means, $\phi(x)$ satisfies the
second equation (Poisson equation) in system $(NSP)$.

As we all know, the knowledge of the solutions for the elliptic system
$$
\left\{ \begin{array}{ll}
  - \Delta u + V(x)u+ \phi u = f(u)  \mbox{ in } \,\,\, \mathbb{R}^{3},\\
-\Delta \phi=u^2  \mbox{ in } \,\,\, \mathbb{R}^{3}, \\
\end{array}
\right.\eqno{(SP)}
$$
 has a great importance in studying stationary solutions
 $\psi(x,t) = e^{-it}u(x)$ of $(NSP)$. It is convenient to observe that the system $(SP)$ contains two kinds of nonlinearities. The first one
is $\phi(x)u$ which is nonlocal, since the electrostatic potential $\phi(x)$ depends also on the wave function,  and is used to describe the interaction between the solitary wave and the electric field. The second type of nonlinearity $f(u)$ is a local one which has been used to model the external forces involving only functions of fields. For more information about the physical background of system $(SP)$, we cite the papers of Benci-Fortunato \cite{bf},
Bokanowski \& Mauser \cite{BM}, Mauser \cite{Mauser}, Ruiz
\cite{Ruiz}, Ambrosetti-Ruiz \cite{AM} and S'anchez \& Soler
\cite{Sanchez}.

An important fact for system $(SP)$ is that it can be reduced into one single Schr\"{o}dinger equation with a nonlocal term (see, for instance, \cite{ap,G,Ruiz,zz}). Effectively, by the Lax-Milgram Theorem, given $u\in H^{1}(\mathbb{R}^{3})$, there exists an unique $\phi=\phi_u \in D^{1,2}(\mathbb{R}^{3})$ such that
$$
-\Delta \phi=u^2 \,\,\, \mbox{in} \,\,\, \mathbb{R}^{3}.
$$
By using standard arguments, we know that $\phi_u$ verifies the following properties (for a proof see
\cite{Daprile1,Ruiz, zz}):
\begin{lemma}{\label{lm1}} For any $u\in H^1(\mathbb{R}^{3})$, we have
\begin{itemize}
\item[$i)$] $\phi_u(x)=\displaystyle \int_{\mathbb{R}^{3}}\frac{u^{2}(y)}{|x-y|}\,dy$ for all $x \in \mathbb{R}^{3}$;
\item[$ii)$]  there exists $C>0$
$$
\int_{\mathbb{R}^{3}}|\nabla \phi_u|^2dx=\int_{\mathbb{R}^{3}}\phi_u u^2dx\leq
C ||u||_{H^{1}(\mathbb{R}^{3})}^4 \quad \forall\, u\in H^1(\mathbb{R}^{3}),
$$
where $||u||_{H^{1}(\mathbb{R}^{3})}=\left(\displaystyle \int_{\mathbb{R}^{3}}(|\nabla u|^2+|u|^{2})dx\right)^{\frac{1}{2}}$ .

\item[$iii)$] $\phi_u\geq 0$ $\forall u\in H^1(\mathbb{R}^{3})$;

\item[$iv)$] $\phi_{tu}=t^2\phi_u$, $\forall t>0$ and $u\in H^1(\mathbb{R}^{3})$;

\item[$v)$] if $u_n \rightharpoonup u$ in $H^1(\mathbb{R}^{3})$, then
$\phi_{u_n} \rightharpoonup \phi_u$ in $D^{1,2}(\mathbb{R}^{3})$ and
$$
\lim_{n \rightarrow +\infty} \int_{\mathbb{R}^{3}}\phi_{ u_n}  u_n ^2dx \geq
\int_{\mathbb{R}^{3}}\phi_u u^2dx.
$$
\item[$vi)$] if $u_n \to  u$ in $H^1(\mathbb{R}^{3})$, then
$\phi_{u_n}\to  \phi_u$ in $D^{1,2}(\mathbb{R}^{3})$. Hence,
$$
\lim_{n \rightarrow +\infty} \int_{\mathbb{R}^{3}}\phi_{ u_n}  u_n ^2dx =
\int_{\mathbb{R}^{3}}\phi_u u^2dx.
$$
\end{itemize}
\end{lemma}

Therefore, $(u,\phi) \in H^1(\mathbb{R}^{3})\times D^{1,2}(\mathbb{R}^{3})$ is a solution of $(SP)$ if, and
only if, $u \in H^{1}(\mathbb{R}^{3})$ is a solution of the nonlocal problem
$$
\left\{ \begin{array}{l}
  - \Delta u + V(x)u+ \phi_u u = f(u) \,\,\, \mbox{ in } \,\,\, \mathbb{R}^{3}, \\
u \in H^{1}(\mathbb{R}^{3}),
\end{array}
\right.\eqno{(P)}
$$
where $\phi_u=\phi \in D^{1,2}(\mathbb{R}^{3})$.

Now, we would like to emphasize that the existence of solutions for problem $(P)$ can be established via variational methods. Associated to the elliptic equation $(P)$, we have the energy functional $I:E \to \mathbb{R}$ given by
$$
I(u)=\frac 12 ||u||^2 +\frac 14\int_{\mathbb{R}^{3}}\phi_u u^2dx-
\int_{\mathbb{R}^{3}}F(u)dx,
$$
where $F(s)=\int_0^s f(t)dt$  and $E$ is the function space
$$
E=\left\{u \in H^{1}(\mathbb{R}^{3})\,:\, \int_{\mathbb{R}^{3}}V(x)|u|^{2}<+\infty \right\}
$$
endowed with the norm
$$
\|u\|=\left(\int_{\mathbb{R}^{3}}(|\nabla u|^{2}+V(x)|u|^{2})dx\right)^{\frac{1}{2}}.
$$

Supposing some conditions on $f$, Lemma \ref{lm1} implies that the functional $I$ is well defined and $I \in C^1(E,\mathbb{R})$ with
\[
I'(u)v =  \int_{\mathbb{R}^{3}}\nabla u \nabla v dx + \int_{\mathbb{R}^{3}}V(x)uv dx+ \int_{\mathbb{R}^{3}}\phi_u
uvdx - \int_{\mathbb{R}^{3}}f(u)vdx \,\,\,\, \forall u,v \in E.
\]
Hence, the critical points of functional $I$ are in fact the weak solutions for nonlocal problem $(P)$.

From  the above commentaries, we know that the system $(SP)$ has a
nontrivial solution if, and only if, the nonlocal problem $(P)$ has a nontrivial
solution. In the last years, many authors had studied the
system $(SP)$ and focused their attentions to establish existence and nonexistence of solutions, multiplicity
of solutions, ground state solutions, radial and nonradial
solutions, semiclassical limit and concentrations of solutions, see for example the papers of
Azzollini \& Pomponio \cite{ap},  Cerami \& Vaira \cite{Cerami},
Coclite \cite{C}, D'Aprile \& Mugnai \cite{Daprile1,Daprile2},
d'Avenia \cite{Davenia}, Ianni \cite{Ianni}, Kikuchi \cite{K1}, and
Zhao \& Zhao \cite{zz}. For the problem set in a bounded
domain, we would like to cite the papers of Siciliano \cite{G}, Ruiz
\& Siciliano \cite{Ruiz-Sic} and Pisani \& Siciliano \cite{LG} for nonnegative solutions and Alves \& Souto \cite{nodal1}, Ianni \cite{Ianni2} and Kim \& Seok \cite{Kim} for sign-changing solutions. However, related to the existence of multi-bump solutions for Schr\"odinger-Poisson system with potential well, as far as we know, there seems to be no existing results.

In the present paper, we will assume that potential $V(x)$ is of the form
$$
V(x)=\lambda a(x)+1,
$$
where $\lambda$ is a positive parameter and $a:\mathbb{R}^{3} \to \mathbb{R}$ is a nonnegative continuous function. Hence, the problems $(SP)$ and $(P)$ can be written respectively as
$$
\left\{ \begin{array}{ll}
  - \Delta u + (\lambda a(x)+1)u+ \phi u = f(u)  \mbox{ in } \,\,\, \mathbb{R}^{3},\\
-\Delta \phi=u^2  \mbox{ in } \,\,\, \mathbb{R}^{3}, \\
\end{array}
\right.\eqno{(SP)_{\lambda}}
$$
and
$$
\left\{ \begin{array}{l}
  - \Delta u + (\lambda a(x)+1)u+ \phi_u u = f(u) \,\,\, \mbox{ in } \,\,\, \mathbb{R}^{3}, \\
u \in H^{1}(\mathbb{R}^{3}).
\end{array}
\right.\eqno{(P)_{\lambda}}
$$

To state the main result, we assume that the function $a(x)$ verifying the following conditions:

\vspace{0,3 cm}

\noindent $(a_1)$ The set $int (a^{-1}(\{0\}))$ is nonempty and there are disjoint open components $\Omega_1, \Omega_2, ....., \Omega_k$ such that
\begin{equation} \label{a1}
int (a^{-1}(\{0\}))= \cup_{j=1}^{k} \Omega_j
\end{equation}
and
\begin{equation} \label{a2}
dist( \Omega_i, \Omega_j)>0 \,\,\, \mbox{for} \,\,\, i \not= j,\,\,\, i, j=1,2,\cdots, k .
\end{equation}

From $(a_1)$, we see that
\begin{equation} \label{E1}
a^{-1}(\{0\})=\cup_{j=1}^{k} \overline{\Omega_j}.
\end{equation}

Related to the function $f$, we will assume the ensuing conditions:
\begin{enumerate}
\item[$(f_1)$\ ]
$\displaystyle\lim_{s\rightarrow 0 } \frac{f(s)}{s} = 0$,
\item[$(f_2)$\ ] $\displaystyle\lim_{|s|\rightarrow +\infty } \frac{f(s)}{s^5} = 0$,
\item[$(f_3)$\ ] There exists $\theta >4$ such that
\[
0 < \theta F(s) \leq sf(s)\quad \forall  s\in \mathbb{R}\setminus \{0\}.
\]
\end{enumerate}
Moreover, we also assume that the nonlinearity $f$ satisfies
\begin{itemize}
\item[$(f_4)$\ ] $\displaystyle \frac{f(s)}{s^3} $ is increasing in $|s|>0$.
\end{itemize}

The motivation to investigate problem $ (SP)_\lambda$ goes back to the papers \cite{Alves} and \cite{DingTanaka}. In \cite{DingTanaka},  inspired by \cite{DelPinoFelmer} and \cite{Sere}, the authors considered the existence of positive multi-bump solution for the problem
\begin{equation} \label{LL2}
\left\{ \begin{array}{ll}
  - \Delta u + (\lambda a(x)+Z(x))u = u^{q} \,\,\,  \mbox{ in } \,\,\, \mathbb{R}^{N},\\
u \in H^{1}(\mathbb{R}^{N}), \\
\end{array}
\right.
\end{equation}
$ q \in (1, \frac{N+2}{N-2} ) $ if $ N \ge 3 $; $ q \in (1, \infty) $ if $ N = 1, 2 $. The authors showed that the above problem has at least $ 2^k-1 $ solutions $u_\lambda$ for large values of $ \lambda $. More precisely, for each non-empty subset $ \Upsilon $ of $ \{ 1,\ldots,k \} $, it was proved that, for any sequence $ \lambda_n \to \infty $ we can extract a subsequence $( \lambda_{n_i}) $ such that $( u_{ \lambda_{n_i} } )$ converges strongly in $ H^1 \big( \mathbb R^N \big) $ to a function $ u $, which satisfies $ u = 0 $ outside $ \Omega_\Upsilon = \bigcup_{ j \in \Upsilon } \Omega_j $ and $ u_{|_{\Omega_j}}, \, j \in \Upsilon $, is a least energy solution for
\begin{equation} \label{LL}
   \begin{cases}
		  - \Delta u + Z(x) u  = u^q, \text{ in } \Omega_j, \\
		  u \in H^1_0 \big( \Omega_j \big), \, u > 0, \text{ in } \Omega_j.
	 \end{cases}
\end{equation}
After, in \cite{Alves}, Alves extended the results described above to the quasilinear  Schr\"odinger equation driven by $p$-Laplacian operator.

Involving the Schr\"odinger-Poisson system with potential wells, there are not so many existing papers. As far as we know, the only paper that considered the existence of solutions for system $(SP)_{\lambda}$ is due to Jiang and Zhou \cite{JZ} where the authors studied the existence and properties of the solutions depending on some parameters. However, nothing is known for the existence of multi-bump type solutions. Motivated by the above references, we intend in the present paper to study the existence of positive multi-bump solution for $(SP)_{\lambda}$. However, we need to point out some difficulties involving this subject:

\vspace{0.3 cm}

\noindent 1- It is well known that the equation \eqref{LL} plays the role of limit equation for \eqref{LL2} as $\lambda$ goes to infinity
and the ground state solution of \eqref{LL} plays an important role in building the multi-bump solutions for \eqref{LL2}. However, little is known about what is the corresponding limit equation for equation $(P)_{\lambda}$ when the parameter $\lambda$ goes to infinity. \\

\noindent 2- Once discovered the limit problem for equation $(P)_{\lambda}$ , it is crucial to prove that it has a specially shaped least energy solution on a subset of the Nehari manifold , see Section 2 for more details. \\

\noindent 3- When we apply variational methods to prove the existence of solution to $(SP)_\lambda$, we are led to study a nonlocal, see problem $(P)_{\lambda}$ above. However, for this class of problem, it is necessary to make a careful revision in the sets used in the deformation lemma found in \cite{Alves} and \cite{DingTanaka} to get multi-bump solution, since they don't work well for this class of system, see Sections 6 and 7 for more details.

\vspace{0.5 cm}

Our main result is the following

\begin{theorem} \label{main}
   Assume that $ (a_1)$ and $ (f_1)-(f_4) $ hold. Then, there exist $ \lambda_0 > 0 $ with the following property: for any non-empty subset $ \Upsilon $ of $ \{1, 2, . . . , k \} $ and $ \lambda \ge \lambda_0 $,  problem $ \big( P_\lambda \big) $ has a positive solution $u_\lambda$. Moreover, if we fix the subset $ \Upsilon $, then for any sequence $ \lambda_n \to \infty $ we can extract a subsequence $ (\lambda_{n_i}) $ such that $ (u_{ \lambda_{n_i} }) $ converges strongly in $ H^{1}(\mathbb R^3) $ to a function $ u $, which satisfies $ u = 0 $ outside $ \Omega_\Upsilon = \cup_{ j \in \Upsilon } \Omega_j $, and $ u_{|_{ \Omega_\Upsilon}}$ is a least energy solution for the nonlocal problem 
$$
\left\{
\begin{array}{l}
 - \Delta u + u+ \left(\displaystyle \int_{\Omega_\Upsilon}\frac{u^2(y)}{|x-y|}dy\right) u = f(u)  \mbox{ in } \,\,\, \Omega_\Upsilon,\\
u(x)>0 \,\,\, \forall x \in \Omega_j \,\,\, \mbox{and} \,\,\, \forall j \in \Upsilon,\\
u \in H^{1}_{0}(\Omega_\Upsilon).
\end{array}
\right.
\eqno{(P)_{\infty, \Upsilon}}
$$
\end{theorem}

In the proof of Theorem \ref{main}, we need to study the existence of least energy solution for problem $(P)_{\infty,\Upsilon}$. The main idea is to prove that the energy function $J$ associated with nonlocal problem $(P)_{\infty, \Upsilon}$ given by
$$
J(u)=\frac 12 \int_{\Omega_\Upsilon}(|\nabla u|^{2}+|u|^{2})dx +\frac 14\int_{\Omega_\Upsilon}\phi_u u^2dx-\int_{\Omega_\Upsilon}F(u)dx,
$$
assumes a minimum value on the set
$$
\mathcal M_{\Upsilon}=\{u\in \mathcal N_{\Upsilon}: J'(u)u_j=0 \mbox { and }
u_{j}\neq 0 \,\,\, \forall j \in \Upsilon \}
$$
where $u_{j}=u_{|_{ \Omega_j}}$ and $\mathcal{N}_{\Upsilon}$ is the corresponding Nehari manifold defined by
$$
\mathcal{N}_{\Upsilon} = \{ u\in
H^1_0(\Omega_\Upsilon)\setminus\{0\}\, :\, J'(u)u=0\}.
$$
More precisely, we will prove that there is $w \in \mathcal{M}_\Upsilon$ such that
$$
J(w)=\inf_{u \in \mathcal{M}_\Upsilon}J(u).
$$
After, we use a deformation lemma to prove that $w$ is a
critical point of $J$, and so, $w$ is a least energy solution
for $(P)_{\infty,\Upsilon}$. The main feature of the least energy solution $w$ is that $w(x)>0 \,\,\, \forall x \in \Omega_j \,\,\, \mbox{and} \,\,\, \forall j \in \Upsilon$ which will be used to describe the existence of multi-bump solutions.

Since we intend to look for positive solutions, through this paper we assume that
$$
f(s)=0, \,\,\, s \leq 0.
$$

\section { The problem $(P)_{\infty,\Upsilon}$ }

In what follows, to show in details the idea of proving the existence of least energy solution for $(P)_{\infty,\Upsilon}$, we will consider $\Upsilon=\{1,2\}$. Moreover, we  will denote by $\Omega$, $\mathcal{N}$ and $ \mathcal M$ the sets $\Omega_\Upsilon$, $\mathcal{N}_\Upsilon$ and $ \mathcal{M}_\Upsilon $ respectively. Thereby,
$$
\Omega=\Omega_1 \cup \Omega_2,
$$
$$
\mathcal{N} = \{ u\in
H^1_0(\Omega)\setminus\{0\}\, :\, J'(u)u=0\}
$$
and
$$
\mathcal M=\{u\in \mathcal N: J'(u)u_1=J'(u)u_2=0 \mbox { and }
u_{1},u_{2}\neq 0 \},
$$
with $u_{j}=u_{|_{ \Omega_j}}$, $j=1,2.$

Since we want to look for least energy for $(P)_{\infty,\Upsilon}$, our goal is to prove the existence of a
critical point for $J$ in the set $\mathcal M$.

\subsection{Technical lemmas}

In what follows, we will denote by $||\,\,\,||$, $||\,\,\,||_1$ and $||\,\,\,||_2$  the norms in $H^{1}_0(\Omega)$, $H^{1}_0(\Omega_1)$ and $H^{1}_0(\Omega_2)$ given by
$$
||u||=\left(\int_{\Omega}(|\nabla u|^{2}+|u|^{2})dx\right)^{\frac{1}{2}},
$$
$$
||u||_1=\left(\int_{\Omega_1}(|\nabla u|^{2}+|u|^{2})dx\right)^{\frac{1}{2}}
$$
and
$$
||u||_2=\left(\int_{\Omega_2}(|\nabla u|^{2}+|u|^{2})dx\right)^{\frac{1}{2}}
$$
respectively.

\vspace{0.5 cm}
In order to show that the set $\mathcal M$ is not empty, we need of the following Lemma.
\begin{lemma}\label{lema3}
Let  $v\in H_0^1(\Omega)$ with $v_{j}\neq 0$ for $j=1,2$. Then,
there are $t,s>0$ such that $J'(tv_1+sv_2)v_1=0$ and
$J'(tv_1+sv_2)v_2=0$.
\end{lemma}

\noindent {\bf{Proof.}} It what follows, we consider the vector field
$$
H(s,t)=\left (J'(tv_1+sv_2)(tv_1),J'(tv_1+sv_2)(sv_2)\right ).
$$
From $(f_1)-(f_3)$, a  straightforward computation yields that there are $0<r<R$ such that
$$
J'(rv_1+sv_2)(rv_1), \,\,\, J'(tv_1+rv_2)(rv_2)>0, \,\,\, \forall s,t\in[r,R]
$$
and
$$
J'(Rv_1+sv_2)(Rv_1), \,\,\, J'(tv_1+Rv_2)(Rv_2)<0, \,\,\, \forall s,t\in[r,R].
$$
Now, the lemma follows by applying  Miranda theorem \cite{Miranda}. \qed

\vspace{0.5 cm}

As an immediate consequence of the last lemma, we have the following corollary

\begin{corollary} \label{CM}

The set $\mathcal M$ is not empty.

\end{corollary}

Next, we will show some technical lemmas.

\begin{lemma}\label{lema2}
There exists $\rho>0$ such that
\begin{enumerate}
\item [(i)] $J(u)\geq ||u||^2/4$ and $||u||\geq \rho, \forall u\in \mathcal N$;
\item [(ii)] $||w_j||_j\geq\rho, \,\, \forall w \in \mathcal M$ and $j=1,2$, where $w_j=w|_{\Omega_j}, j=1,2$.
\end{enumerate}
\end{lemma}

\noindent {\bf{Proof.}} From $(f_4)$, for any $u\in \mathcal N$
$$
4J(u)=4J(u)-J'(u)u= ||u||^2+ \int_{\Omega}[uf(u)-4F(u)]dx\geq
||u||^2
$$
and so,
$$
J(u)\geq ||u||^2/4, \,\,\, \forall u \in \mathcal N .
$$

From $(f_1)$ and $(f_2)$, there is $C>0$ such that
$$
f(s)s\leq \frac {\lambda_1}2 s^2+Cs^6, \mbox { for all } s\in
\mathbb R,
$$
where $\lambda_1$ is the first eigenvalue of $(-\Delta , H_{0}^{1}(\Omega)).$ Since $J'(u)u=0$,
$$
||u||^2<||u||^2+\int_\Omega \phi_{u}u^2dx=\int_{\Omega}uf(u)dx\leq
\frac {\lambda_1}2\int_{\Omega}u^2dx +C \int_{\Omega}u^6dx.
$$
Then, by Sobolev embeddings,
$$
||u||^2< \frac{1}{2} ||u||^2+\hat C ||u||^6,
$$
from where it follows that
$$
 ||u||\geq \rho  \,\,\, \forall u \in \mathcal N,
$$
where $\rho = \left(\frac 1{2\hat C}\right)^{\frac{1}{4}}$, finishing the proof of $(i).$

If $w \in \mathcal M$, we have that $J'(w)w_1=J'(w)w_2=0$. Then, a simple computation gives $J'(w_{j})w_{j}<0$ for $j=1,2$, which implies
$$
||w_{j}||_j^2<||w_{j}||_j^2+\int_{\Omega_j} \phi_{w_{j}}(w_{j})^2dx<\int_{\Omega_j}
f(w_j)w_jdx, \,\,\, \mbox{for} \,\,\, j=1,2.
$$
As in $(i)$, we can deduce that $||w_{j}||_j\geq \rho$ for $j=1,2$.

 \qed

 \begin{lemma}\label{lema2x}
If $(w_n)$ is a bounded sequence in  $\mathcal M$ and $p\in
(2,6)$, we have
$$
\liminf_n \int_{\Omega_j}|w_{n,j}|^{p}dx>0 \,\,\, j=1,2.
$$
where $w_{n,j}=w_n|_{\Omega_j}$ for $j=1,2$. \end{lemma}

\noindent {\bf{Proof.}} From $(f_1)$ and $(f_2)$, given
$\varepsilon>0$ there exists $C>0$ such that
$$
f(s)s\leq \varepsilon\lambda_1 s^2+C|s|^p+\varepsilon s^6, \mbox {
for all } s\in \mathbb R.
$$
Since $w_n\in \mathcal M $, by Lemma \ref{lema2}
$$
\rho^{2} \leq||w_{n,j}||_j^2<\int_{\Omega_j}w_{n,j}f(w_{n,j})dx \leq
\varepsilon\lambda_1\int_{\Omega_j}(w_{n,j})^2dx
+C\int_{\Omega_j}|w_{n,j}|^pdx+\varepsilon \int_{\Omega_j}(w_{n,j})^6dx
$$
that is,
$$
\rho^2 \leq \varepsilon \left( \lambda_1\int_{\Omega_j}(w_{n,j})^2dx + \int_{\Omega_j}(w_{n,j})^6dx  \right) + C\int_{\Omega_j}|w_{n,j}|^pdx.
$$
Using the boundedness of $(w_n)$, there is $C_1$ such that
$$
\rho^2 \leq \varepsilon C_1 + C\int_{\Omega_j}|w_{n,j}|^pdx.
$$
Fixing $\varepsilon = \frac{\rho^2}{2C_1}$, we get
$$
\int_{\Omega_j}|w_{n,j}|^pdx \geq \frac{\rho^2}{2C},
$$
showing that
$$
\liminf_{n} \int_{\Omega_j}|w_{n,j}|^pdx \geq \frac{\rho^2}{2C} >0.
$$
\qed

\subsection { Existence of least energy solution for $(P)_{\infty,\Upsilon}$  }

In this subsection, our main goal is to prove the following result

\begin{theorem} \label{T2}
Assume that $(f_1)-(f_4)$ hold. Then equation $(P)_{\infty,\Upsilon}$ possesses a positive least energy solution on the set $\mathcal{M}$.
\end{theorem}

\noindent {\bf{Proof.}} In what follows, we denote by $c_0$ the infimum of $J$ on $\mathcal M$, that
is,
$$
c_0=\inf_{v\in \mathcal M} J(v).
$$
From Lemma \ref{lema2}(i), we conclude that $c_0>0$.

By Corollary \ref{CM}, we know that $\mathcal M$ is not empty, then there is a sequence $(w_n) \subset \mathcal M$ satisfying
$$
\lim_{n} J(w_n)=c_0.
$$
Still from Lemma \ref{lema2}(i), $(w_n)$ is a bounded
sequence. Hence, without loss of generality, we may suppose that there is $w \in H_0^1(\Omega)$ verifying
$$
w_n \rightharpoonup w \,\,\, \mbox{in} \,\,\, H_0^{1}(\Omega),
$$

$$
w_n \to w \,\,\, \mbox{in} \,\,\,  L^p(\Omega) \,\,\ \forall \, p \in [1,2^{*})
$$
and
$$
w_n(x) \to w(x) \,\,\, \mbox{a.e. in} \,\, \Omega.
$$
Then, $(f_2)$ combined with the \emph{compactness lemma of Strauss}
\cite[Theorem A.I, p.338]{bl} gives
$$
\lim_n \int_{\Omega_j}|w_{n,j}|^{p}dx = \int_{\Omega_j}|w_j|^{p}dx,
$$
$$
\lim_n \int_{\Omega_j} w_{n,j}f(w_{n,j})dx=\int_{\Omega_j} w_j f(w_j)dx
$$
and
$$
\lim_n \int_{\Omega_j}  F(w_{n,j})dx=\int_{\Omega_j} F(w_j)dx,
$$
from where it follows together with Lemma \ref{lema2x} that $w_j\neq 0$ for $j=1,2$. Then, by Lemma \ref{lema3} there are $t,s>0$ verifying
$$
J'(tw_1+sw_2)w_1=0 \,\,\, \mbox{and} \,\,\, J'(tw_1+sw_2)w_2=0.
$$
Next, we will show that $t,s\leq 1$. Since $J'(w_{n,j})w_{n,j}=0$ for $j=1,2$,
$$
||w_{n,1}||_1^2+\int_{\Omega_1} \phi_{w_{n,1}}(w_{n,1})^2dx+\int_{\Omega_1}
\phi_{w_{n,2}}(w_{n,1})^2dx=\int_{\Omega_1} f(w_{n,1})w_{n,1}dx
$$
and
$$
||w_{n,2}||_2^2+\int_{\Omega_2} \phi_{w_{n,2}}(w_{n,2})^2dx+\int_{\Omega_2}
\phi_{w_{n,1}}(w_{n,2})^2dx=\int_{\Omega_2} f(w_{n,2})w_{n,2}dx.
$$
Taking the limit in the above equalities, we obtain
$$
||w_1||_1^2+\int_{\Omega_1} \phi_{w_1}(w_1)^2dx+\int_{\Omega_1}
\phi_{w_2}(w_1)^2dx\leq\int_\Omega f(w_1)w_1dx
$$
and
$$
||w_2||_2^2+\int_{\Omega_2} \phi_{w_2}(w_2)^2dx+\int_{\Omega_2}
\phi_{w_1}(w_2)^2dx\leq\int_{\Omega_2} f(w_2)w_2dx.
$$
Recalling that
$$
J'(tw_1+sw_2)(tw_1)=J'(tw_1+sw_2)(sw_2)=0,
$$
it follows that
$$
t^2||w_1||_2^2+t^4\int_{\Omega_1} \phi_{w_1}(w_1)^2dx+t^2s^2\int_{\Omega_1}
\phi_{w_2}(w_1)^2dx=\int_{\Omega_1} f(tw_1)tw_1dx
$$
and
$$
s^2||w_2||_2^2+s^4\int_{\Omega_2} \phi_{w_2}(w_2)^2dx+t^2s^2\int_{\Omega_2}
\phi_{w_2}(w_1)^2dx=\int_{\Omega_2} f(sw_2)sw_2dx.
$$
Now, without loss of generality, we will suppose that $s\geq t$.
Under this condition,
$$
s^2||w_2||_2^2+s^4\int_{\Omega_2} \phi_{w_2}(w_2)^2dx+s^4\int_{\Omega_2}
\phi_{w_2}(w_1)^2dx\geq\int_{\Omega_2} f(sw_2)sw_2dx
$$
and then
$$
\left( \frac 1{s^2}-1\right )||w_2||_2^2\geq\int_{\Omega_2} \left ( \frac
{f(sw_2)sw_2}{(sw_2)^4}-\frac {f(w_2)w_2}{(w_2)^4} \right
)(w_2)^4dx.
$$
If $s>1$, the left side in this inequality is negative, but from
$(f_4)$, the right side is positive, thus we must have $s \leq 1$, which also implies that $t \leq 1$.

Our next step is to show that $J(tw_1+sw_2)=c_0$.  Recalling that $tw_1+sw_2\in
\mathcal M$, we derive that
$$
c_0 \leq J(tw_1+sw_2)=J(tw_1+sw_2)-\frac{1}{4}
J'(tw_1+sw_2)(tw_1+sw_2).
$$
Hence,
$$
c_0 \leq \left (J(tw_1)-\frac 14
J'(tw_1)(tw_1)\right )+\left (J(sw_2)-\frac 14
J'(sw_2)(sw_2)\right ).
$$
From  $(f_4)$,
$$
J(tw_1)-\frac 14 J'(tw_1)(tw_1)  \leq J(w_1)-\frac 14 J'(w_1)(w_1)
$$
and
$$
J(sw_2)-\frac 14 J'(sw_2)(sw_2) \leq J(w_2)-\frac 14 J'(w_2)(w_2),
$$
leading to
$$
c_0 \leq \left (J(w_1)-\frac 14
J'(w_1)(w_1)\right )+\left (J(w_2)-\frac 14
J'(w_2)(w_2)\right ).
$$
Using Fatous' Lemma together with $(f_4)$, we see that
$$
c_0 \leq  J(tw_1 + sw_2) \leq \liminf_{n} \left (J(w_n) -\frac{1}{4}J'(w_n)w_n \right)=\lim_n J(w_n)=c_0,
$$
which means that
$$
c_0 = J(tw_1 + sw_2).
$$

Until this moment, we have proved that there exists a $w_o=tw_1 + sw_2\in
\mathcal M$, such that $J(w_o)=c_0$. In what follows, let us denote $w_o$ by $w$, consequently
$$
J(w)=c_0 \,\,\, \mbox{and} \,\,\, w\in \mathcal M.
$$

To complete the proof of Theorem \ref{main}, we claim that $w$ is a critical point for
functional $J$.  To see why, for each $\varphi \in H^{1}_{0}(\Omega)$, we introduce the functions $Q^{i}:\mathbb{R}^{3} \to \mathbb{R},i=1,2$ given by
$$
\begin{array}{l}
Q^{1}(r,z,l)=\displaystyle \int_{\Omega_1}|\nabla(w_1+r\varphi_1+zw_1)|^{2}dx+\int_{\Omega_1}\phi_{(w+r\varphi+zw_1+lw_2)}(w_1+r\varphi_1+zw_1)^{2}dx\\
\mbox{}\\
\hspace{2 cm} -\displaystyle \int_{\Omega_1}f(w_1+r\varphi_1+zw_1)(w_1+r\varphi_1+zw_1 )dx
\end{array}
$$
and
$$
\begin{array}{l}
Q^{2}(r,z,l)=\displaystyle \int_{\Omega_2}|\nabla(w_2+r\varphi_2+lw_2)|^{2}dx+\int_{\Omega_2}\phi_{(w+r\varphi+zw_1+lw_2)}(w_2+r\varphi_2+lw_2)^{2}dx\\
\mbox{}\\
\hspace{2 cm} -\displaystyle \int_{\Omega_2}f(w_2+r\varphi_2+lw_2)(w_2+r\varphi_2+lw_2 )dx.
\end{array}
$$
By a direct computation,
$$
\frac{\partial Q^{1}}{\partial z}(0,0,0)=2\int_{\Omega_1}|\nabla w_1|^{2}dx+4\int_{\Omega_1}\phi_w w_1^{2}dx-\int_{\Omega_1}(f'(w_1)w_1^{2}+f(w_1)w_1)  dx
$$
and so,
$$
\frac{\partial Q^{1}}{\partial z}(0,0,0)= \int_{\Omega_1}(f(w_1)w_1-f'(w_1)w_1^{2})dx+2\int_{\Omega_1}\phi_w w_1^{2}dx.
$$
By $(f_4)$, we know that $f'(s)s^{2} \geq 3 f(s)s$  for all $ s\geq 0$.  Thus,
$$
\frac{\partial Q^{1}}{\partial z}(0,0,0) \leq -2\left(\int_{\Omega_1}f(w_1)w_1dx-\int_{\Omega_1}\phi_w w_1^{2}dx\right).
$$
Now, recalling that $J'(w)w_1=0$, we have
$$
\|w_1\|^{2}+ \int_{\Omega_1}\phi_w w_1^{2}dx=\int_{\Omega_1}f(w_1)w_1dx.
$$
Then,
$$
\frac{\partial Q^{1}}{\partial z}(0,0,0) \leq -2\|w_1\|^{2}<0.
$$
The same type of argument gives
$$
\frac{\partial Q^{2}}{\partial l}(0,0,0) \leq -2\|w_2\|^{2}<0.
$$
Therefore, applying the implicit function theorem, there are functions $z(r), l(r)$ of class $C^{1}$ defined on some interval $(-\delta, \delta), \delta>0$ such that $z(0)=l(0)=0$ and
$$
Q^{i}(r,z(r),l(r))=0, \,\,\, r \in (-\delta, \delta), i=1,2.
$$
This shows that for any $t \in (-\delta, \delta)$,
$$
v(t)=w+r\varphi +z(r)w_1 +l(r)w_2 \in {\mathcal M}.
$$
Since,
$$
J(w)=c_0=\inf_{v \in {\mathcal M}}J(v),
$$
we derive that
$$
J(v(r)) \geq J(w), \,\,\,\, \forall r \in (-\delta, \delta),
$$
that is,
$$
J(w+r\varphi +z(r)w_1 +l(r)w_2) \geq J(w), \,\,\,\, \forall r \in (-\delta, \delta).
$$
From this,
$$
\frac{J(w+r\varphi +z(r)w_1 +l(r)w_2) - J(w) }{r} \geq 0, \,\,\,\, \forall r \in (0, \delta).
$$
Taking the limit of $r \to 0$, we get
$$
J'(w)(\varphi+z'(0)w_1+l'(0)w_2) \geq 0.
$$
Recalling that $J'(w)w_1=J'(w)w_2=0$, the above inequality loads to
$$
J'(w)\varphi \geq 0, \,\,\, \forall \varphi \in H^{1}_{0}(\Omega)
$$
and so,
$$
J'(w)\varphi= 0, \,\,\, \forall \varphi \in H^{1}_{0}(\Omega),
$$
showing that $w$ is a critical point for $J$. \qed

\section{An auxiliary problem}

In this section, we work with an auxiliary problem adapting the ideas explored by del Pino \& Felmer in \cite{DelPinoFelmer} (see also \cite{Alves} and \cite{DingTanaka}).

We start recalling that the energy functional $ I_\lambda \colon E_\lambda \to \mathbb R $ associated with $ (P)_\lambda $ is  given by
$$
   I_\lambda (u) = \frac{1}{2}\int_{ \mathbb R^3 }  \left( | \nabla u |^{2 } + \big( \lambda a(x) + 1 \big) u^{ 2 } \right)dx + \frac{1}{4}\int_{\mathbb{R}^{3}}\phi_u u^{2}dx- \int_{ \mathbb R^3 } F(u)dx,
$$
where $ E_\lambda = \big( E, \| \cdot \|_\lambda \big) $ with
$$
   E = \left\{ u \in H^{1} ( \mathbb R^3 ) \, ; \, \int_{ \mathbb R^3 } a(x) |u|^{ 2 }dx < \infty \right\},
$$
and
$$
   \| u \|_\lambda = \left(\int_{\mathbb{R}^{3}}(|\nabla u|^{2}+(\lambda a(x)+1)|u|^{2})dx \right)^{\frac{1}{2}}.
$$
Thus $ E_\lambda \hookrightarrow H^{1}( \mathbb R^3 ) $ continuously for $ \lambda \geq 0 $ and $ E_\lambda $ is compactly embedded in $ L_{ loc }^{s}( \mathbb R^3 ) $, for all $ 1 \leq s \leq 2^*=\frac{2N}{N-2} $ for $N \geq 3$. A  direct computation gives that $ E_\lambda $ is a Hilbert space. Also, being $ {\cal O} \subset \mathbb R^3 $ an open set, from the relation
\begin{equation} \label{modular relation 1}
  \int_{ \cal O } \left( \big| \nabla u \big|^{ 2} + \big( \lambda a(x) + 1 \big) | u |^{ 2 } \right)dx  \ge \int_{ \cal O } |u|^{ 2 }dx ,
\end{equation}
for all $ u \in E_\lambda $ with $ \lambda \geq 0 $, fixed $ \delta \in (0,1)$, there is $ \nu > 0 $, such that
\begin{equation} \label{modular relation 2}
   \|u\|_{ \lambda, \cal{O} }^{2} - \nu |u|_{2, \cal{O} }^{2} \geq \delta \|u\|_{ \lambda,\cal{O} }^{2}, \, \forall u \in E_\lambda, \, \lambda \geq 0.
\end{equation}
Hereafter,
$$
\|u\|_{ \lambda, \cal{O} }= \left(\int_{ \cal O }( | \nabla u |^{ 2} + ( \lambda a(x) + 1 ) | u |^{ 2 })dx \right)^{\frac{1}{2}}
$$
and
$$
 |u|_{2, \cal{O} }= \left(\int_{ \cal O } | u |^{ 2 }dx \right)^{\frac{1}{2}}.
$$

We recall that for any $ \epsilon > 0 $, the hypotheses $ (f_1) $ and $ (f_2) $ yield
\begin{equation} \label{f estimate}
   f(s) \le \epsilon | s | + C_\epsilon | s|^{2^{*}-1}, \, \forall x \in \mathbb R^3 \,\,\, \mbox{and} \,\,\, s \in \mathbb R.
\end{equation}
Consequently,
\begin{equation} \label{F estimate}
   F(s) \le \epsilon | s |^{2 } + C_\epsilon |s|^{2^{*}}, \, \forall x \in \mathbb R^3 \,\,\, \mbox{and} \,\,\, s  \in \mathbb R,
\end{equation}
where $ C_\epsilon $ depends on $ \epsilon $. Moreover, for  $\nu >0$ fixed in (\ref{modular relation 2}), the assumptions $ (f_1) $ and $ (f_4) $ imply that there is an unique $a>0$ verifying
\begin{equation} \label{NU}
\frac{f(a)}{a}=\nu
\end{equation}

Using the numbers $a$ and $\nu$, we set the function $ \tilde{f} \colon \mathbb R^3 \times \mathbb R \to \mathbb R $ given by		 
$$
   \tilde{f}(s) =
   \begin{cases}
      \ \, f(s), \ s \le a \\
	    \nu \, s, \ s \ge a
   \end{cases},
$$
which fulfills the inequality
\begin{equation} \label{til f estimate}
   \tilde{f}(s) \le \nu | s|, \,\,\, \forall  s \in \mathbb R.
\end{equation}
Thus
\begin{equation} \label{t til f estimate}
   \tilde{f}(s) s \le \nu | s |^{2}, \,\,\, \forall s \in \mathbb R
\end{equation}
and
\begin{equation} \label{til F estimate}
  \tilde{F}(s) \le \frac{\nu}{2} | t |^{ 2 }, \, \forall s \in \mathbb R,
\end{equation}
where $ \tilde{F}(s) = \int_0^t \tilde{f}(t) \, dt $.

Now, since $ \Omega=int( a^{ -1 } (\{0\})) $ is formed by $ k $ connected components $ \Omega_1, \ldots, \Omega_k $ with $ \text{dist} \big( \Omega_i, \Omega_j \big) > 0, \, i \ne j $, then for each $ j \in \{ 1, \ldots, k \} $, we are able to fix a smooth bounded domain $ \Omega'_j $ such that
\begin{equation} \label{omega}
   \overline{\Omega_j} \subset \Omega'_j \, \text{ and } \, \overline{\Omega'_i} \cap \overline{\Omega'_j} = \emptyset, \text{ for } i \ne j.
\end{equation}

From now on, we fix a non-empty subset $ \Upsilon \subset \left\{ 1, \ldots, k \right\} $ and
$$
   \Omega_\Upsilon = \bigcup_{ j \in \Upsilon } \Omega_j, \, \Omega'_\Upsilon = \bigcup_{ j \in \Upsilon } \Omega'_j \,\, \mbox{and} \,\,\,
	 \chi_\Upsilon =
	 \begin{cases}
	    1, \text{ if } x \in \Omega'_\Upsilon \\
		  0, \text{ if } x \notin \Omega'_\Upsilon .
	 \end{cases}
$$
Using the above notations, we set the functions
$$
   g(x,s) = \chi_\Upsilon(x) f(s) + \big( 1-\chi_\Upsilon(x) \big) \tilde{f}(s), \, (x,s) \in \mathbb R^3 \times \mathbb R
$$
and
$$
   G(x,s) = \int_0^s g(x,t) \, dt, \, (x,s) \in \mathbb R^3 \times \mathbb R,
$$
and the auxiliary nonlocal problem
$$
\left\{
\begin{array}{l}
- \Delta u + ( \lambda a(x) + 1 ) u + \phi_u u = g(x,u), \text{ in } \mathbb R^3, \\
u \in E_\lambda .
\end{array}
\right.
\leqno{(A_\lambda)}
$$

The problem $ \big( A_\lambda \big) $ is related to $ \big( P_\lambda \big) $ in the sense that, if $ u_\lambda $ is a solution for $ \big( A_\lambda \big) $  verifying
$$
    u_\lambda (x) \leq a, \, \forall x \in \mathbb R^N \setminus \Omega'_\Upsilon,
$$
then it is a solution for $ \big( P_\lambda \big) $.

In comparison to $ \big( P_\lambda \big) $, problem $ \big( A_\lambda \big) $ has the advantage that the energy functional associated with $ \big( A_\lambda \big) $, namely, $ \phi_\lambda \colon E_\lambda \to \mathbb R $ given by
$$
    \phi_\lambda(u) = \frac{1}{2} \int_{ \mathbb R^3 } ( \left| \nabla u \right|^{ 2 } + ( \lambda a(x) + 1 ) | u |^{ 2 } )dx +\frac{1}{4}\int_{\mathbb{R}^{3}}\phi_u u^{2}dx- \int_{ \mathbb R^3} G(x,u)dx,
$$
satisfies the $ (PS) $ condition, whereas $ I_\lambda $ does not necessarily satisfy this condition.

\vspace{0.5 cm}

\begin{proposition} \label{boundedness}
   All $ (PS)_d $ sequences for $ \phi_\lambda $ are bounded in $ E_\lambda $.
\end{proposition}

\noindent {\bf Proof.}
   Let $ (u_n) $ be a $ (PS)_d $ sequence for $ \phi_\lambda $. So, there is $ n_0 \in \mathbb N $ such that
$$
   \phi_\lambda (u_n) - \frac{1}{\theta} \phi_\lambda'(u_n) u_n \le d+1 + \| u_n \|_\lambda, \text { for } n \ge n_0.
$$
On the other hand, by (\ref{t til f estimate}) and (\ref{til F estimate})
$$
   \tilde{F}(s) - \frac{1}{\theta} \tilde{f}(s)s \le \left( \frac{1}{2} - \frac{1}{\theta} \right) \nu | s |^{ 2 }, \, \forall x \in \mathbb R^3, s \in \mathbb R,
$$
which together with (\ref{modular relation 2}) gives
$$
   \phi_\lambda (u_n) - \frac{1}{\theta} \phi_\lambda'(u_n) u_n \ge \left( \frac{1}{2} - \frac{1}{\theta} \right) \delta \|u_n\|^{2}_{\lambda}, \, \forall n \in \mathbb N,
$$
from where it follows that $ (u_n) $ is bounded in $ E_\lambda $.
\qed

\begin{proposition} \label{Estimativa no infinito}
If $(u_n)$ is a $(PS)_d$ sequence for $\phi_{\lambda}$, then given $\epsilon>0$, there is $R>0$ such that
\begin{equation} \label{Estimativa}
   \limsup_n \int_{ \mathbb R^3 \setminus B_R (0) } (| \nabla u_n|^{2 } + ( \lambda a(x) + 1 ) | u_n |^{ 2 } )dx < \epsilon.
\end{equation}
Hence, once that $g$ has a subcritical growth, if $ u \in E_\lambda $ is the weak limit of $ (u_n) $, then
$$
\int_{\mathbb R^3}g(x,u_n)u_n\,dx \to \int_{\mathbb R^3} g(x,u)u \, dx \, \text{ and } \, \int_{\mathbb R^3} g(x,u_n)v \, dx \to \int_{\mathbb R^3} g(x,u)v \, dx, \, \forall v \in E_\lambda.
$$
\end{proposition}

\noindent {\bf Proof.}
   Let $ (u_n) $ be a $ (PS)_d $ sequence for $ \phi_\lambda $, $ R > 0 $ large such that $ \Omega'_\Upsilon \subset B_{ \frac{R}{2} }(0) $ and $ \eta_R \in C^\infty (\Bbb R^3) $ satisfying
$$
   \eta_R (x) =
	 \begin{cases}
	    0, \, x \in B_{ \frac{R}{2} }(0) \\
			1, \, x \in \Bbb R^3 \setminus B_R (0)
	 \end{cases},
$$
$ 0 \le \eta_R \le 1 $ and $ \big| \nabla \eta_R \big| \le \dfrac{C}{R} $, where $ C > 0 $ does not depend on $ R $. This way,
\begin{align*}
   \mbox{}  & \int_{ \Bbb R^3 } (  | \nabla u_n \big|^{2 } + ( \lambda a(x) + 1 ) | u_n |^{ 2 } ) \eta_R dx+ \int_{\mathbb{R}^{3}}\phi_{u_n}u_n^{2}\eta_Rdx\\
	   = & \phi_\lambda'(u_n) \left( u_n \eta_R \right) - \int_{ \Bbb R^3 } u_n  \nabla u_n \nabla \eta_R dx+ \int_{ \Bbb R^3 \setminus \Omega'_\Upsilon } \tilde{f}(u_n) u_n \eta_Rdx.
\end{align*}
Denoting
$$
L =\int_{ \Bbb R^3 } (  | \nabla u_n |^{ 2 } + ( \lambda a(x) +1 ) | u_n |^{2 } ) \eta_R dx,
$$
it follows from (\ref{t til f estimate}),
$$
   L \leq \phi_\lambda'(u_n) ( u_n \eta_R ) + \frac{C}{R} \int_{ \Bbb R^3 } | u_n | | \nabla u_n |dx + \nu \int_{ \Bbb R^3 } | u_n |^{2 } \eta_R dx.
$$
Using H\"older's inequality, we derive
$$
	 L \le \phi_\lambda'(u_n) \left( u_n \eta_R \right) + \frac{C}{R} | u_n |_{ 2} | \nabla u_n |_2+ \nu L.
$$	
Since $ (u_n) $ and $ ( | \nabla u_n | ) $ are bounded in $ L^{2}( \Bbb R^3 ) $, we obtain
$$
L \leq o_n(1) + \frac{C}{(1-\nu)R}.
$$
Therefore
$$
   \limsup_n \int_{ \Bbb R^3 \setminus B_R(0) } (  | \nabla u_n |^{2 } + ( \lambda a(x) +1 ) | u_n |^{2} )dx \le \frac{C}{(1-\nu)R}.
$$
So, given $ \epsilon > 0 $,  choosing a $ R > 0 $ possibly still bigger, we have that $ \dfrac{C}{(1-\nu)R} < \epsilon $, which proves (\ref{Estimativa}). Now, we will show that
$$
\int_{\Bbb R^3}g(x,u_n)u_ndx \to \int_{\Bbb R^3}g(x,u)udx.
$$
Using the fact that $g(x,u)u \in L^{1}(\Bbb R^3)$ together with (\ref{Estimativa}) and Sobolev embeddings, given  $\epsilon >0$, we can choose $R>0$ such that
$$
\limsup_{n \to +\infty}\int_{\Bbb R^3 \setminus B_R(0)}|g(x,u_n)u_n|dx \leq \frac{\epsilon}{4} \quad \mbox{and} \quad \int_{\Bbb R^3 \setminus B_R(0)}|g(x,u)u|dx \leq \frac{\epsilon}{4}.
$$
On the other hand, since $g$ has a subcritical growth, we have by compact embeddings
$$
\int_{B_{R}(0)}g(x,u_n)u_ndx \to \int_{B_{R}(0)}g(x,u)udx.
$$
Combining the above information, we conclude that
$$
\int_{\Bbb R^3}g(x,u_n)u_ndx \to \int_{\Bbb R^3}g(x,u)udx.
$$
The same type of arguments works to prove that
$$
\int_{\Bbb R^3}g(x,u_n)v dx\to \int_{\Bbb R^3}g(x,u)vdx, \quad \forall v \in E_{\lambda}.
$$

\qed

\begin{proposition} \label{(PS) condition}
   $ \phi_\lambda $ verifies the $ (PS) $ condition.
\end{proposition}
		
\noindent {\bf Proof}
Let $ (u_n) $ be a $ (PS)_d $ sequence for $ \phi_\lambda $ and $ u \in E_\lambda $ such that $u_n \rightharpoonup u$ in $E_{\lambda}$. Thereby, by Proposition \ref{Estimativa no infinito},
$$
   \int_{ \mathbb R^3 } g(x,u_n)u_ndx \to \int_{ \mathbb R^3 } g(x,u)u dx\, \text{ and } \, \int_{ \mathbb R^3 } g(x,u_n)v dx\to \int_{ \mathbb R^3 } g(x,u)vdx, \, \forall v \in E_\lambda.
$$
Moreover, the weak limit also gives
$$
\int_{ \mathbb R^3 }  \nabla u \nabla ( u_n-u )dx  \to 0
$$
and
$$
\int_{ \mathbb R^3 }  ( \lambda a(x) + 1 ) u ( u_n-u )dx  \to 0.
$$
Recalling that $\phi_\lambda'(u_n)u_n=o_n(1)$  and $\phi_\lambda'(u_n)u=o_n(1)$, the above limits lead to
$$
 \|u_n-u\|^{2}_\lambda \to 0,
$$
finishing the proof. \qed

%------------------------------------------------------------THE (PS) INFINITY CONDITION-------------------------------------------------------------------------------------------

\section{The $ (PS)_\infty $ condition}

%-----------------------------------------------------------------------------------------------------------------------------------------------------------------------------------

A sequence $ (u_n) \subset  H^{ 1} ( \mathbb R^3 ) $ is called a $ (PS)_\infty $ \emph{sequence for the family} $ \left( \phi_\lambda \right)_{\lambda \ge 1} $, if there is a sequence $ ( \lambda_n ) \subset  [1, \infty) $ with $ \lambda_n \to \infty $, as $ n \to \infty $, verifying
$$
   \phi_{ \lambda_n }(u_n) \to c \text{ and } \left\| \phi'_{ \lambda_n }(u_n) \right\|_{E^{*}_{\lambda_n}} \to 0, \text{ as } n \to \infty,
$$
for some $c \in \mathbb{R}$.
\begin{proposition} \label{(PS) infty condition}
   Let $ (u_n) \subset H^{ 1, } ( \mathbb R^3 ) $ be a $ (PS)_\infty $ sequence for $ \left( \phi_\lambda \right)_{\lambda \ge 1} $. Then, up to a subsequence, there exists $ u \in H^{ 1 } ( \mathbb R^3 ) $ such that $ u_n \rightharpoonup u $ in $ H^{ 1}( \mathbb R^3 ) $. Furthermore,
\begin{enumerate}
  \item[(i)] $ u_n \to u $ in $ H^{ 1} ( \mathbb R^3 ) $;
	\item[(ii)] $ u = 0 $ in $ \mathbb R^3 \setminus \Omega_\Upsilon $,  $ u_{|_{\Omega_j}} \geq 0$ for all $j \in \Upsilon $, and $u$ is a solution for
   $$
   	    \begin{cases}
		     - \Delta u + u + \phi_u u = f(u), \text{ in } \Omega_\Upsilon, \\
				   u \in H^{1}_0 ( \Omega_\Upsilon );
			\end{cases}
   \eqno{(P)_{\infty,\Upsilon} }
	$$
	 \item[(iii)] $ \displaystyle \lambda_n\int_{\mathbb R^3}  a(x) | u_n |^{ 2 } \to 0 $;
	 \item[(iv)] $ \|u_n-u\|^{2}_{\lambda,\Omega'_\Upsilon}\to 0, \text{ for } j \in \Upsilon $;
   \item[(v)] $  \|u_n\|^{2}_{\lambda,\mathbb{R}^{3} \setminus \Omega'_\Upsilon} \to 0 $;
	 \item[(vi)] $ \phi_{ \lambda_n } (u_n) \to \displaystyle \frac{1}{2}\int_{ \Omega_\Upsilon } (| \nabla u |^{ 2 } +  | u |^{2} )dx +\frac{1}{4}\int_{\Omega_{\Upsilon}}\phi_{u} u^{2}dx- \int_{ \Omega_\Upsilon } F(u)dx $.
\end{enumerate}
\end{proposition}

\noindent {\bf Proof.}
   Using the Proposition \ref{boundedness}, we know that $( \| u_n \|_{ \lambda_n } ) $ is bounded in $ \mathbb R $ and $ (u_n) $ is bounded in $ H^{1}( \mathbb R^3 ) $. So, up to a subsequence, there exists $ u \in H^{1}(\mathbb R^3) $ such that
$$
   u_n \rightharpoonup u  \text{ in } H^{1} ( \mathbb R^3) \, \text{ and } \, u_n(x) \to u(x) \text{ for a.e. } x \in \mathbb R^3.
$$
Now, for each $ m \in \mathbb N $, we define $ C_m = \left\{ x \in \mathbb R^3 \, ; \, a(x) \ge \dfrac{1}{m} \right\} $. Without loss of generality, we can assume $ \lambda_n < 2 ( \lambda_n-1 ), \, \forall n \in \mathbb N $. Thus
	       $$
		       \int_{ C_m } | u_n |^{ 2 }dx \le \frac{2m}{\lambda_n} \int_{ C_m } \big( \lambda_n a(x)+1) | u_n |^{ 2 }dx \leq \frac{C}{\lambda_n}.
	       $$
By Fatou's lemma, we derive
	       $$
			     \int_{ C_m } | u |^{ 2 }dx = 0,
	       $$
	       which implies that $ u = 0 $ in $ C_m $ and, consequently, $ u = 0 $ in $ \mathbb R^3 \setminus \overline{\Omega} $. From this, we are able to prove $(i)-(vi)$.
	
\begin{enumerate}
   \item[$(i)$] Since $ u = 0 $  in $ \mathbb R^3 \setminus \overline{\Omega} $, repeating the argument explored in  Proposition \ref{(PS) condition} we get
	       $$
            \int_{ \mathbb R^3 } (|\nabla u_n - \nabla u|^{2}  + ( \lambda_n a(x) + 1) |u_n-u|^{2})dx \to 0,
         $$
  which implies $ u_n \to u $ in $ H^{1}(\mathbb R^3) $.
	 \item[$(ii)$] Since $ u \in H^{ 1}(\mathbb R^3) $ and $ u = 0 $  in $ \mathbb R^3 \setminus \overline{\Omega} $, we have $ u \in H^{1}_0( \Omega ) $ or, equivalently, $ u_{ |_{\Omega_j} } \in H^{1}_0( \Omega_j) $, for $ j = 1, \ldots, k $. Moreover, the limit $u_n \to u$ in $H^{1}(\mathbb R^3)$ combined with $\phi'_{\lambda_n}(u_n)\varphi \to 0$  for $\varphi \in C^{\infty}_0 ( \Omega_\Upsilon )$ implies that
         \begin{equation} \label{u is solution}	
	          \int_{\Omega_\Upsilon}( \nabla u  \nabla \varphi +  u \varphi )dx +\int_{\Omega_\Upsilon}\phi_u \varphi dx - \int_{\Omega_\Upsilon} f(u) \varphi dx = 0,
         \end{equation}
showing that $ u_{ |_{\Omega_{\Upsilon}} }  $ is a solution for the nonlocal problem
         $$
				 \begin{cases}
		        - \Delta u + u +\phi_u u= f(u), \text{ in } \Omega_\Upsilon, \\
		        u \in H^{1}_0 ( \Omega_\Upsilon ).
		       \end{cases}
					\eqno{(P)_{\infty,\Upsilon} }
$$
	
 On the other hand, if $ j \notin \Upsilon $, we must have
				 $$
				   \int_{ \Omega_j }(|\nabla u|^{2} + |u|^{2}) dx+\int_{ \Omega_j } \phi_u u^{2} dx- \int_{ \Omega_j } \tilde{f}(u)u dx= 0.
				 $$
				The above equality combined with (\ref{t til f estimate}) and (\ref{modular relation 2}) gives
			 	 $$
				    0 \ge \|u\|^{2}_{ \lambda, \Omega_j }- \nu \|u\|^{2}_{2,\Omega_j } \geq \delta \|u\|^{2}_{ \lambda, \Omega_j }(u) \geq 0,
				 $$
				 from where it follows $ u_{|_{\Omega_j}} = 0 $ for $ j \notin \Upsilon $.  This proves $ u = 0 $ outside $ \Omega_\Upsilon $ and $ u \geq 0 $ in $ \Bbb R^3$.
				
	\item[$(iii)$]  From (i),
	       $$
		        \lambda_n \int_{ \mathbb R^3 }  a(x) | u_n |^{2 } dx= \int_{ \mathbb R^2 } \lambda_n a(x) | u_n-u |^{ 2 }dx \leq \|u_n-u\|^{2}_{\lambda_n},
	       $$
				implying that
				$$
				\lambda_n \int_{ \mathbb R^3 }  a(x) | u_n |^{2 } dx \to 0.
				$$
	\item[$(iv)$] Let $ j \in \Upsilon $. From (i),
	       $$
				    |u_n-u|^{2}_{2, \Omega'_j }, |\nabla u_n - \nabla u|^{2}_{2, \Omega'_j } \to 0.
				 $$
				 Then,
				 $$
				 \int_{ \Omega'_\Upsilon } ( | \nabla u_n |^{ 2 } - | \nabla u |^{ 2 } )dx \to 0 \quad \mbox{and} \quad \int_{ \Omega'_\Upsilon } ( | u_n |^{ 2 } - | u |^{ 2 } )dx \to 0.
				 $$
				 From (iii),
				 $$
				    \int_{ \Omega'_\Upsilon } \lambda_n a(x) | u_n |^{2}  dx \to 0.
			   $$
				 This way
				 $$
				  \|u_n\|^{2}_{ \lambda_n, \Omega'_\Upsilon } \to \int_{ \Omega_\Upsilon } ( | \nabla u |^{ 2 } + | u |^{ 2} )dx.
				 $$
	 \item[$(v)$] By (i),  $ \|u_n-u\|^{2}_{ \lambda_n } \to 0 $, and so,
	              $$
								  \|u_n\|^{2}_{ \lambda_n, \mathbb R^3 \setminus \Omega_\Upsilon } \to 0.
								$$
	 \item[$(vi)$] We can write the functional  $\phi_{\lambda_n}$ in the following way
	       \begin{multline*}
	          \phi_{ \lambda_n } (u_n) = \sum_{ j \in \Upsilon } \left[\frac{1}{2} \int_{ \Omega'_j }(| \nabla u_n |^{2} + ( \lambda_n a(x) + 1) | u_n |^{2})dx +\frac{1}{4}\int_{\Omega'_j}\phi_{u_n} u_n^{2}dx\right]\\
						+ \frac{1}{2} \int_{ \mathbb R^3\setminus \Omega'_\Upsilon } ( | \nabla u_n |^{ 2 } + ( \lambda_n a(x) +1 ) | u_n |^{ 2 } ) dx+\frac{1}{4}\int_{\mathbb R^3\setminus \Omega'_\Upsilon}\phi_{u_n} u_n^{2} dx- \int_{ \mathbb R^3 } G(x,u_n)dx.
				 \end{multline*}
				 From $(i)-(v)$,
				 $$
					  \frac{1}{2}\int_{ \Omega'_j } ( | \nabla u_n |^{ 2 } + ( \lambda_n a(x) + 1 ) | u_n |^{2 } )dx \to \frac{1}{2}\int_{ \Omega_j }(| \nabla u|^{2} + | u |^{2})dx,
				 $$

				 $$
					 \frac{1}{2}\int_{ \mathbb R^3 \setminus \Omega'_\Upsilon } (| \nabla u_n |^{2} + ( \lambda_n a(x) + 1) | u_n |^{2} ) dx\to 0,
				$$
				$$
			  \int_{\Omega'_j}\phi_{u_n} u_n^{2}dx \to \int_{\Omega_j}\phi_{u} u^{2}dx,
				$$
				$$
				\int_{\mathbb R^3\setminus \Omega'_\Upsilon}\phi_{u_n} u_n^{2} \to 0
				$$
				and
				 $$
					  \int_{ \mathbb R^3 } G(x,u_n)dx \to \int_{ \Omega_\Upsilon } F(u)dx.
				 $$
				 Therefore
				 $$
					\phi_{ \lambda_n } (u_n) \to \frac{1}{2}\int_{ \Omega_\Upsilon } ( | \nabla u |^{2} + | u |^{2} )dx+\frac{1}{4}\int_{\Omega_{\Upsilon}}\phi_{u} u^{2} dx- \int_{ \Omega_\Upsilon } F(u)dx.
				 $$
\end{enumerate}	
\qed

%----------------------------------------------------------THE BOUNDEDNESS OF THE $ (A_\LAMBDA ) $ SOLUTIONS-----------------------------------------------------------------------

\section{The boundedness of the $ \big( A_\lambda \big) $ solutions}

%-----------------------------------------------------------------------------------------------------------------------------------------------------------------------------------

In this section, we study the boundedness outside $ \Omega'_\Upsilon $ for some solutions of $ \big( A_\lambda \big ) $. To this end, we adapt the arguments found in \cite{Alves} and \cite{Gon} for our new setting.

\begin{proposition} \label{P:boundedness of the solutions}
    Let $ \big( u_\lambda \big) $ be a family of solutions for $ \big( A_\lambda \big) $ such that $ u_\lambda \to 0 $ in $ H^{1}( \mathbb R^3 \setminus \Omega_\Upsilon ) $, as $ \lambda \to \infty $. Then, there exists $ \lambda^* > 0 $ with the following property:
$$
   \left| u_\lambda \right|_{ \infty, \mathbb R^3 \setminus \Omega'_\Upsilon } \le a, \, \forall \lambda \ge \lambda^*.
$$
Hence, $u_{\lambda}$ is a solution for $(P_\lambda)$ for $\lambda \geq \lambda^*$.
\end{proposition}

\vspace{.3cm}
\noindent {\bf Proof.}
Since $\partial \Omega'_\Upsilon$ is a compact set, fixed a neighborhood $\mathcal B$ of $\partial \Omega'_\Upsilon$ such that
$$
\mathcal B \subset \mathbb{R}^{3} \setminus \Omega_\Upsilon,
$$
the interation Moser technique implies that there is $C>0$, which is independent of $\lambda$, such that
$$
|u_\lambda|_{L^{\infty}(\partial \Omega'_\Upsilon) }
 \leq C |u_{\lambda}|_{L^{2^{*}}({\mathcal B})}
$$
Since $u_\lambda \to 0$ in $H^{1}(\mathbb{R}^{3} \setminus \Omega_\Upsilon)$, we have that
$$
|u_{\lambda}|_{L^{2^{*}}({\mathcal B})} \to 0 \,\,\, \mbox{as} \,\,\, \lambda \to \infty.
$$
Hence, there is $\lambda^{*}>0$ such that
$$
|u_{\lambda}|_{L^{2^{*}}({\mathcal B})} < \frac{a}{2C} \,\,\, \forall \lambda \geq \lambda^{*},
$$
and so,
$$
|u_\lambda|_{L^{\infty}(\partial \Omega'_\Upsilon) }< a \,\,\, \forall \lambda \geq \lambda^{*}.
$$
Next, for $ \lambda \ge \lambda^* $, we set  $ \widetilde{u}_\lambda \colon \mathbb R^3 \setminus \Omega'_\Upsilon \to \mathbb R $ given by
$$
   \widetilde{u}_\lambda(x) = ( u_\lambda-a )^+ (x).
$$
Thereby,  $ \widetilde{u}_\lambda \in H^{1}_0(\mathbb R^3 \setminus \Omega'_\Upsilon) $. Our goal is showing that  $\widetilde{u}_\lambda = 0 $ in $ \mathbb R^3 \setminus \Omega'_\Upsilon $. This implies
$$
   \left| u_\lambda \right|_{ \infty, \mathbb R^3 \setminus \Omega'_\Upsilon } \leq a.
$$
In fact, extending $ \widetilde{u}_\lambda = 0 $ in $ \Omega'_\Upsilon $ and taking $ \widetilde{u}_\lambda $ as a test function, we obtain
$$
   \int_{ \mathbb R^3 \setminus \Omega'_\Upsilon } \nabla u_\lambda  \nabla \widetilde{u}_\lambda dx + \int_{ \mathbb R^3 \setminus \Omega'_\Upsilon } \! \! \! \! ( \lambda a(x) + 1)u_\lambda \widetilde{u}_\lambda dx  \leq \int_{ \mathbb R^N \setminus \Omega'_\Upsilon } g \left( x, u_\lambda \right) \widetilde{u}_\lambda dx.
$$
Since
\begin{gather*}
   \int_{ \mathbb R^3 \setminus \Omega'_\Upsilon } \nabla u_\lambda \nabla \widetilde{u}_\lambda dx = \int_{ \mathbb R^3 \setminus \Omega'_\Upsilon }| \nabla \widetilde{u}_\lambda |^{2} dx, \\
	\int_{ \mathbb R^3 \setminus \Omega'_\Upsilon } ( \lambda a(x) +1) u_\lambda \widetilde{u}_\lambda dx= \int_{ \left( \mathbb R^3 \setminus \Omega'_\Upsilon \right)_+ } ( \lambda a(x) + 1) \left( \widetilde{u}_\lambda+a \right) \widetilde{u}_\lambda dx
\end{gather*}
and
$$
  \int_{ \mathbb R^3 \setminus \Omega'_\Upsilon } g \left( x, u_\lambda \right) \widetilde{u}_\lambda dx = \int_{ \left( \mathbb R^3 \setminus \Omega'_\Upsilon \right)_+ } \frac{g \left( x, u_\lambda \right)}{u_\lambda} \left( \widetilde{u}_\lambda+a \right) \widetilde{u}_\lambda dx,
$$
where
$$
   \left( \mathbb R^3 \setminus \Omega'_\Upsilon \right)_+ = \left\{ x \in \mathbb R^3 \setminus \Omega'_\Upsilon \, ; \, u_\lambda(x) > a \right\},
$$
we derive
$$
   \int_{ \mathbb R^3 \setminus \Omega'_\Upsilon } | \nabla \widetilde{u}_\lambda|^{2 }dx + \int_{ \left( \mathbb R^3 \setminus \Omega'_\Upsilon \right)_+ } \! \! \! \! ( (\lambda a(x) + 1 ) -\frac{g \left( x, u_\lambda \right)}{u_\lambda}) \left( \widetilde{u}_\lambda+a \right)  \widetilde{u}_\lambda dx \leq 0,
$$
Now, by \eqref{til f estimate},
$$
 ( \lambda a(x) + 1) - \frac{g ( x, u_\lambda )}{u_\lambda} > \nu  - \frac{\tilde{f} \left( x, u_\lambda \right)}{u_\lambda} \ge 0 \quad \mbox{in} \quad  \left( \mathbb R^3 \setminus \Omega'_\Upsilon \right)_+ .
$$
This form, $ \widetilde{u}_\lambda = 0 $ in $ \left( \mathbb R^3 \setminus \Omega'_\Upsilon \right)_+ $. Obviously, $ \widetilde{u}_\lambda = 0  $ at the points where $ u_\lambda \leq a $, consequently, $ \widetilde{u}_\lambda = 0 $ in $ \mathbb R^N \setminus \Omega'_\Upsilon $. 		
\qed

%--------------------------------------------------A SPECIAL CRITICAL VALUE FOR PHI_\LAMBDA-----------------------------------------------------------------------------------------

\section{A special minimax value for $ \phi_\lambda $}

%-----------------------------------------------------------------------------------------------------------------------------------------------------------------------------------

For fixed non-empty subset $ \Upsilon \subset \left\{ 1, \ldots, k \right\} $ , consider
$$
   I_\Upsilon(u) = \frac{1}{2}\int_{\Omega_\Upsilon} (| \nabla u |^{2} + | u |^{2} ) dx+ \frac{1}{4}\int_{\Omega_\Upsilon}\phi_{u} u^{2}dx  - \int_{ \Omega_\Upsilon } F(u)dx, \ u \in H^{1}_0( \Omega_\Upsilon), 
$$
the energy functional associated to $ (P)_{\infty,\Upsilon} $,  and $\phi_{ \lambda,\Upsilon }: H^{1}(\Omega'_\Upsilon) \to \mathbb{R} $ given by 
$$
   \phi_{ \lambda,\Upsilon }(u) = \frac{1}{2}\int_{ \Omega'_\Upsilon } (| \nabla u |^{2} +( \lambda a(x) + 1) | u |^{2} )dx +  \frac{1}{4}\int_{\Omega'_\Upsilon}\left(\displaystyle \int_{\Omega'_\Upsilon}\frac{\tilde{u}^2}{|x-y|}dy\right) u^{2}dx - \int_{ \Omega'_\Upsilon } F(u)dx, 
$$
the energy functional associated to the nonlocal problem 
$$
   \begin{cases}
	    - \Delta u + ( \lambda a(x) + 1) u + \left(\displaystyle \int_{\Omega'_\Upsilon}\frac{\tilde{u}^2}{|x-y|}dy\right) u = f(u), \text{ in } \Omega'_\Upsilon, \\
			\frac{\partial u}{\partial \eta} = 0, \text{ on } \partial \Omega'_\Upsilon,
	 \end{cases}
$$
where
$$
\tilde{u}(x)=
\left\{
\begin{array}{l}
u(x), \,\,\,\, \text{in} \,\,\, \Omega'_{\Upsilon} \\
0, \,\,\,\, \text{in} \,\,\, \mathbb{R}^{3} \setminus \Omega'_\Upsilon .
\end{array}
\right.
$$

In the following, we denote by $c_\Upsilon$ the number given by
$$
c_\Upsilon=\inf_{u \in \mathcal M_{\Upsilon}}I_\Upsilon(u)
$$
where
$$
\mathcal M_{\Upsilon}=\{u\in \mathcal N_{\Upsilon}: I_\Upsilon'(u)u_j=0 \mbox { and }
u_{j}\neq 0 \,\,\, \forall j \in \Upsilon \}
$$
with $u_{j}=u_{|_{ \Omega_j}}$ and
$$
\mathcal{N}_{\Upsilon} = \{ u\in
H^1_0(\Omega_\Upsilon)\setminus\{0\}\, :\, I_\Upsilon'(u)u=0\}.
$$
Of a similar way,  we denote by $c_{\lambda,\Upsilon}$ the number given by
$$
c_{\lambda,\Upsilon}=\inf_{u \in \mathcal M'_{\Upsilon}}\phi_{\lambda,\Upsilon}(u)
$$
where
$$
\mathcal M'_{\Upsilon}=\{u\in \mathcal N'_{\Upsilon}: \phi_{\lambda,\Upsilon}'(u)u_j=0 \mbox { and }
u_{j}\neq 0 \,\,\, \forall j \in \Upsilon \}
$$
with $u_{j}=u_{|_{ \Omega'_\Upsilon}}$ and
$$
\mathcal{N}'_{\Upsilon} = \{ u\in
H^1(\Omega'_\Upsilon)\setminus\{0\}\, :\,\phi_{\lambda,\Upsilon}'(u)u=0\}.
$$
Repeating the same approach used in Section 2,  we ensure that there exist $ w_\Upsilon  \in H^{1}_0( \Omega_\Upsilon  ) $ and $ w_{ \lambda,\Upsilon  } \in H^{1}( \Omega'_\Upsilon ) $ such that
$$
   I_\Upsilon( w_\Upsilon) = c_\Upsilon  \, \text{ and } \, I'_\Upsilon ( w_\Upsilon) = 0
$$
and
$$
	\phi_{ \lambda,\Upsilon}( w_{ \lambda,\Upsilon }) = c_{ \lambda,\Upsilon}	\, \text{ and } \, \phi'_{ \lambda,\Upsilon }( w_{ \lambda,\Upsilon }) = 0.
$$

By a direct computation, it is possible to show that there is $\tau>0$ such that if $u \in \mathcal M_{\Upsilon}$, then
\begin{equation} \label{tau0}
\|u_j\|_j > \tau, \,\,\, \forall j \in \Upsilon,
\end{equation}
where, $\|\,\,\,\,\|_j$ denotes the norm on $H^{1}_0(\Omega_j)$ given by
$$
||u||_j=\left(\int_{\Omega_j}(|\nabla u|^{2}+|u|^{2})dx\right)^{\frac{1}{2}}.
$$
In particular, since $w_\Upsilon \in \mathcal M_{\Upsilon}$, we also have
\begin{equation} \label{tau1}
\|w_{\Upsilon,j}\|_j > \tau \,\,\, \forall j \in \Upsilon,
\end{equation}
where $w_{\Upsilon,j}=w_{\Upsilon}|_{\Omega_j}$ for all $j \in \Upsilon$.
\vspace{.3cm}
\begin{lemma}
   There holds that
\begin{enumerate}
   \item[(i)] $ 0 < c_{ \lambda,\Upsilon  } \le c_\Upsilon , \, \forall \lambda \geq 0 $;
   \item[(ii)] $ c_{ \lambda,\Upsilon  } \to c_\Upsilon , \text{ as } \lambda \to \infty, $.
\end{enumerate}
\end{lemma}

\noindent {\bf Proof}
\begin{enumerate}
   \item[(i)]  Using the inclusion $ H^{1}_0(\Omega_\Upsilon) \subset H^{1}(\Omega'_\Upsilon) $, it  easy to observe that
				 $$
				   c_{ \lambda,\Upsilon  } \leq  c_\Upsilon .
				 $$
	 \item[(ii)] It suffices to show that $ c_{ \lambda_n,j } \to c_j, \text{ as } n \to \infty $, for all sequences $ ( \lambda_n ) $ in $ [1,\infty) $ with $ \lambda_n \to \infty, \text{ as } n \to \infty  $. Let $ \left( \lambda_n \right) $ be such a sequence and consider an arbitrary subsequence of $ \left( c_{ \lambda_n,\Upsilon  } \right) $ (not relabelled) . Let $ w_n \in H^{1}( \Omega'_j) $ with % = w_{ \lambda_n,j }
		     $$
					  \phi_{ \lambda_n,\Upsilon  }(w_n) = c_{ \lambda_n,\Upsilon }	\, \text{ and } \, \phi'_{ \lambda_n,\Upsilon } ( w_n) = 0.
				 $$
By the previous item, $ \big( c_{ \lambda_n,\Upsilon  } \big) $ is bounded. Then, there exists $( w_{ n_k }) $ subsequence of $( w_n) $ such that $ (\phi_{ \lambda_{ n_k },\Upsilon }( w_{ n_k }) )$ converges and  $\phi'_{ \lambda_{ n_k },\Upsilon } ( w_{ n_k }) = 0 $. Now, repeating similar arguments explored in the proof of  Proposition \ref{(PS) infty condition}, there is $ w \in H^{1}_0(\Omega_\Upsilon ) \setminus \{0\} \subset H^{1}(\Omega'_\Upsilon ) $ such that
$$
w_j =w|_{\Omega_j} \not=0  \,\,\, \forall j \in \Upsilon
$$	
and 				
				
				$$
				    w_{ n_k } \to w \text{ in } H^{1}( \Omega'_\Upsilon ), \text{ as } k \to \infty.
				 $$
				 Furthermore, we also can prove that
				 $$
				    c_{ \lambda_{ n_k },\Upsilon } = \phi_{ \lambda_{ n_k },\Upsilon  }( w_{ n_k }) \to I_\Upsilon (w)
				 $$
				and
				 $$
				    0 = \phi'_{ \lambda_{ n_k },\Upsilon  }( w_{ n_k }) \to I'_\Upsilon (w).
				 $$
				Then,  $ w \in \mathcal M_{\Upsilon}$, and  by definition of $c_{\Upsilon}$,
				 $$
				    \lim_k c_{ \lambda_{ n_k },\Upsilon  } \geq c_\Upsilon .
				 $$
				 The last inequality together with item (i) implies
				 $$
				    c_{ \lambda_{ n_k },\Upsilon  } \to c_\Upsilon , \text{ as } k \to \infty.
				 $$
				 This establishes the asserted result.
\end{enumerate}
\qed
\vspace{.5cm}

In the sequel,  we fix  $ R > 1 $ verifying
\begin{equation} \label{R}
   0< I'_j \left( \frac{1}{R} w_j \right)\left( \frac{1}{R} w_j \right) \,\,\, \mbox{and} \,\,\, I'_j(R w_j)(R w_j)<0,  \text{ for } j \in \Upsilon,
\end{equation}
where $I_j$ denotes the ensuing functional
$$
I_{ j}(u) = \frac{1}{2}\int_{ \Omega_j } (| \nabla u |^{2} +| u |^{2} )dx +\frac{1}{4}\int_{\Omega_j}\phi_u u^{2}dx - \int_{ \Omega_j } F(u)dx, \ u \in H^{1}_0(\Omega_j)
$$
with $\phi_u$ being the solution of the problem
$$
\left\{
\begin{array}{l}
-\Delta \phi =u^{2}, \,\,\, \text{in} \,\,\, \mathbb{R}^{3} \\
\phi \in D^{1,2}(\mathbb{R}^{3}).
\end{array}
\right.
$$

In the sequel, to simplify the notation, we rename the components $ \Omega_j $ of $ \Omega $ in way such that $ \Upsilon = \{ 1, 2, \ldots, l \} $  for some $ 1 \le l \le k $. Then, we define:
\begin{gather*}
   \gamma_0 ( \textbf{t} )(x) = \sum_{j=1}^l t_j R w_j(x) \in H^{1}_0(\Omega_\Upsilon), \, \forall \textbf{t}=( t_1, \ldots, t_l )\in [1/R^2,1]^l, \\
   \Gamma_\ast = \Big\{ \gamma \in C \big( [1/R^2,1]^l, E_\lambda \setminus \{ 0 \} \big) \, ; \, \gamma(\textbf{t})|_{\Omega'_j} \not=0 \,\,\, \forall j \in \Upsilon\, ; \; \gamma = \gamma_0 \text{ on } \partial [1/R^2,1]^l \Big\}
\end{gather*}
and
$$
	 b_{ \lambda, \Upsilon } = \inf_{ \gamma \in \Gamma_\ast } \max_{ \textbf{t}\in [1/R^2,1]^l } \phi_\lambda \big( \gamma ( \textbf{t} ) \big).
$$

\vspace{.3cm}
Next, our intention is proving an important relation among $ b_{ \lambda, \Upsilon } $, $c_\Upsilon$ and $c_{\lambda, \Upsilon}$. However, to do this, we need to some technical lemmas. The arguments used are the same found in \cite{Alves}, however for reader's convenience we will repeat their proofs

\begin{lemma} \label{solution's existence}
   For all $ \gamma \in \Gamma_\ast $, there exists $ (s_1, \ldots, s_l ) \in [1/R^2,1]^l $ such that
$$
   \phi'_{ \lambda,j } \big( \gamma ( s_1, \ldots, s_l ) \big) \big( \gamma ( s_1, \ldots, s_l ) \big) = 0, \, \forall j \in \Upsilon
$$
where
$$
\phi_{ \lambda,j}(u) = \frac{1}{2}\int_{ \Omega'_j } (| \nabla u |^{2} +( \lambda a(x) + 1) | u |^{2} )  dx +\frac{1}{4}\int_{\Omega'_j}\phi_{\tilde u} u^{2} dx- \int_{ \Omega'_\Upsilon } F(u)dx, \ u \in H^{1}(\Omega'_\Upsilon).
$$
\end{lemma}

\noindent {\bf Proof}
   Given $ \gamma \in \Gamma_\ast $, consider $ \widetilde{\gamma} \colon [1/R^2,1]^l \to \mathbb R^l $ such that
$$
   \widetilde{\gamma} ( \textbf{t} ) = \Big( \phi'_{ \lambda,1 } \big( \gamma ( \textbf{t} ) \big) \gamma ( \textbf{t} ), \ldots, \phi'_{ \lambda,l } \big( \gamma ( \textbf{t} ) \big) \gamma ( \textbf{t} ) \Big), \text{ where } \textbf{t} = ( t_1, \ldots, t_l ).
$$
For $ \textbf{t} \in \partial [1/R^2,1]^l $, it holds $ \widetilde{\gamma} ( \textbf{t} ) = \widetilde{\gamma_0} ( \textbf{t} ) $. From this, we observe that there is no $ \textbf{t} \in \partial [1/R^2,1]^l  $ with $ \widetilde{\gamma} ( \textbf{t} ) = 0 $. Indeed, for any $ j \in \Upsilon $,
$$
   \phi'_{ \lambda,j} ( \gamma_0 ( {\bf{t}} )) \gamma_0 ( {\bf{t}} ) = I'_j( \gamma_0 ( {\bf{t}} )) \gamma_0 ( {\bf{t}} ).
$$

This form, if $ {\bf{t}} \in \partial [1/R^2,1]^l $, then $ t_{j_0} =1 $ or $ t_{j_0} = \frac{1}{R^2} $, for some $ j_0 \in \Upsilon $. Consequently,
$$
   \phi'_{ \lambda,j_0 } \big( \gamma_0 ( {\bf{t}} ) \big) \gamma_0 ( {\bf{t}} ) = I'_{j} ( R w_{j_0} ) ( R w_{j_0} ) <0
$$
or
$$
\phi'_{ \lambda,j_0 } \big( \gamma_0 ( {\bf{t}} ) \big) \gamma_0 ( {\bf{t}} ) = I'_{j} \left( \frac{1}{R} w_{j_0} \right) \left( \frac{1}{R} w_{j_0} \right)>0.
$$
Therefore, from \eqref{R},
$$
\phi'_{ \lambda,j_0 } \big( \gamma_0 ( {\bf{t}} ) \big) \gamma_0 ( {\bf{t}} ) \not= 0.
$$
Now, we able to compute the degree $ \deg \big( \widetilde{\gamma}, (1/R^2,1)^l, (0, \ldots, 0 ) \big) $. Since
$$
   \deg \big( \widetilde{\gamma}, (1/R^2,1)^l, (0, \ldots, 0 ) \big) = \deg \big( \widetilde{\gamma_0}, (1/R^2,1)^l, (0, \ldots, 0 ) \big),
$$
and, for $ \textbf{t} \in (1/R^2,1)^l $,
$$
  \widetilde{\gamma_0} ( \textbf{t} ) = 0 \iff {\bf{t}} = \left( \frac{1}{R}, \ldots, \frac{1}{R} \right),
$$
we derive
$$
   \deg \big( \widetilde{\gamma}, (1/R^2,1)^l, (0, \ldots, 0 ) \big) \ne 0.
$$
This shows what was stated.
\qed

\begin{proposition} \label{blambdagamma}
  \mbox{}
\begin{enumerate}
   \item[(i)] $ c_{ \lambda,\Upsilon } \le b_{ \lambda,\Upsilon } \le c_\Upsilon, \, \forall \lambda \ge 1 $;
	 \item[(ii)] $ b_{ \lambda,\Upsilon } \to c_\Upsilon, \text{ as } \lambda \to \infty $;
	 \item[(iii)] $ \phi_\lambda \big( \gamma({\bf{t}}) \big) < c_\Upsilon, \, \forall \lambda \ge 1, \gamma \in \Gamma_\ast \text{ and } {\bf{t}} = (t_1, \ldots, t_l ) \in \partial [1/R^2,1]^l $.
\end{enumerate}	
\end{proposition}

\noindent {\bf Proof}
\begin{enumerate}
   \item[(i)] Since $ \gamma_0 \in \Gamma_\ast $,
	       $$
				    b_{ \lambda,\Upsilon } \le \max_{ ( t_1, \ldots, t_l ) \in [1/R^2,1]^l } \phi_\lambda ( \gamma_0 ( t_1, \ldots, t_l ) ) \leq  \max_{ ( t_1, \ldots, t_l )\in \mathbb{R}^l } I_\Upsilon ( \sum_{j=1}^{l}t_j R w_j ) = c_\Upsilon.
				 $$
				 Now, fixing $ {\bf s} = (s_1, \ldots, s_l) \in [1/R^2,1]^l $ given in Lemma \ref{solution's existence} and recalling that
				 $$
				    c_{ \lambda,\Upsilon } =\inf_{u \in \mathcal M'_{\Upsilon}}\phi_{\lambda,\Upsilon}(u)
$$
where
$$
\mathcal M'_{\Upsilon}=\{u\in \mathcal N'_{\Upsilon}: \phi_{\lambda,\Upsilon}'(u)u_j=0 \mbox { and }
u_{j}\neq 0 \,\,\, \forall j \in \Upsilon \},
$$
				 it follows that
				 $$
				    \phi_{ \lambda,\Upsilon } ( \gamma( {\bf s } )) \geq c_{ \lambda,\Upsilon }.
				 $$
			From \eqref{til F estimate},
				 $$
				    \phi_{ \lambda, \mathbb R^3 \setminus \Omega'_\Upsilon } (u) \geq 0, \, \forall u \in H^{1} \big( \mathbb R^3 \setminus \Omega'_\Upsilon \big),
				 $$
which leads to
				 $$
				    \phi_\lambda \big( \gamma( {\textbf t} ) \big) \geq \phi_{ \lambda,\Upsilon } ( \gamma( {\textbf t} )) , \, \forall \textbf{t} = (t_1, \ldots, t_l) \in [1/R^2,1]^l.
				 $$
				 Thus
				 $$
				    \max_{ ( t_1, \ldots, t_l )\in [1/R^2,1]^l } \phi_\lambda \big( \gamma ( t_1, \ldots, t_l ) \big) \geq \phi_{\lambda,\Upsilon} \big( \gamma( \textbf{s} ) \big) \ge c_{ \lambda,\Upsilon },
				 $$
			showing that
				 $$
				    b_{ \lambda,\Upsilon } \ge c_{ \lambda,\Upsilon };
				 $$
	 \item[(ii)] This limit is clear by the previous item, since we already know $ c_{ \lambda,\Upsilon } \to c_\Upsilon $, as $ \lambda \to \infty $;
	 \item[(iii)] For $ \textbf{t} = ( t_1, \ldots, t_l ) \in \partial [1/R^2,1]^l $, it holds $ \gamma ( \textbf{t} ) = \gamma_0 ( \textbf{t} ) $. From this,
	       $$
				    \phi_\lambda \big( \gamma ( \textbf{t} ) \big) = I_\Upsilon ( \gamma_0 ( \textbf{t} )).
				 $$
				 From \eqref{R},  we derive
				 $$
				    \phi_\lambda \big( \gamma ( \textbf{t} ) \big) \le c_\Upsilon - \epsilon,
				 $$
				 for some $ \epsilon > 0 $, so (iii) holds.
\end{enumerate}
\qed

%----------------------------------------------------------------THE PROOF OF THE MAIN THEOREM--------------------------------------------------------------------------------------

\section{Proof of the main theorem}

%-----------------------------------------------------------------------------------------------------------------------------------------------------------------------------------

To prove Theorem \ref{main}, we need to find nonnegative solutions $ u_\lambda $ for large values of $ \lambda $, which converges to a least energy solution of $(P)_{\infty,\Upsilon}$ as $ \lambda \to \infty $. To this end, we will show two propositions which together with the Propositions \ref{(PS) infty condition} and  \ref{P:boundedness of the solutions} will imply that Theorem \ref{main} holds.

Henceforth, we denote by
$$
\Theta=\left\{u \in E_\lambda\,:\, \|u\|_{\lambda, \Omega'_j}> \frac{\tau}{8R} \,\,\, \forall j \in \Upsilon \right\}
$$
and
$$
	\phi_\lambda^{ c_\Upsilon } = \big\{ u \in E_\lambda \, ; \, \phi_\lambda(u) \le c_{ \Upsilon } \big\}.
$$

Moreover, we fix $\delta=\frac{\tau}{48R}$, and for $ \mu >0 $,
\begin{equation} \label{A}
   {\cal A}_\mu^\lambda = \left\{ u \in \Theta_{2\delta} \,:\, | \phi_{ \lambda}(u)-c_\Upsilon | \leq \mu\right\}.
\end{equation}
We observe that
$$
   w_\Upsilon \in {\cal A}_\mu^\lambda \cap \phi_\lambda^{ c_\Upsilon },
$$
showing that $ {\cal A}_\mu^\lambda \cap \phi_\lambda^{ c_\Upsilon } \ne \emptyset $.
we have the following uniform estimate of $ \big\| \phi'_{ \lambda }(u) \big\|_{E^{*}_{\lambda}} $ on the region $ \left( {\cal A}_{ 2 \mu }^\lambda \setminus {\cal A}_\mu^\lambda \right) \cap \phi_\lambda^{ c_\Upsilon } $.

\begin{proposition} \label{derivative estimate}
For each  $ \mu > 0 $, there exist $ \Lambda_\ast \ge 1 $ and $ \sigma_0 >0   $ independent of $ \lambda $ such that
\begin{equation}
   \big\| \phi'_{ \lambda }(u) \big\|_{E^{*}_{\lambda}} \ge \sigma_0, \text{ for } \lambda \ge \Lambda_\ast \text{ and all } u \in \left( {\cal A}_{ 2 \mu }^\lambda \setminus {\cal A}_\mu^\lambda \right) \cap \phi_\lambda^{ c_\Upsilon }.
\end{equation}
\end{proposition}

\noindent {\bf Proof}
   We assume that there exist $ \lambda_n \to \infty $ and $ u_n \in \left( {\cal A}_{ 2 \mu }^{\lambda_n} \setminus {\cal A}_\mu^{\lambda_n} \right) \cap \phi_{\lambda_n}^{ c_\Upsilon } $ such that
$$
   \big\| \phi'_{ \lambda_n }(u_n) \big\|_{E^{*}_{\lambda_n}} \to 0.
$$
Since $ u_n \in {\cal A}_{ 2 \mu }^{ \lambda_n } $, this implies $( \|u_n\|_{ \lambda_n }) $ is a bounded sequence and, consequently, it follows that $ \big( \phi_{ \lambda_n }(u_n) \big) $ is also bounded. Thus, passing a subsequence if necessary, we can assume $(\phi_{ \lambda_n }(u_n)) $ converges. Thus, from Proposition \ref{(PS) infty condition}, there exists $ 0 \le u \in H^{1}_0( \Omega_\Upsilon ) $ such that $u$ is a solution for $ (SP)_\Upsilon $,
$$
u_n \to u \,\, \text{in} \,\, H^{1}(\mathbb{R}^{3}), \,\,	 \|u_n\|_{ \lambda_n, \mathbb R^3 \setminus \Omega_\Upsilon } \to 0 \, \text{ and } \,  \phi_{ \lambda_n} (u_n) \to I_\Upsilon(u).
$$
Recalling that $(u_n) \subset \Theta_{2 \delta}$, we derive that
$$
\|u_n\|_{\lambda_n,\Omega'_j} > \frac{\tau}{12R} \,\,\, \forall j \in \Upsilon.
$$
Then, taking the limit of $n \to +\infty$, we find
$$
\|u\|_j \geq \frac{\tau}{12R} \,\,\, \forall j \in \Upsilon,
$$
yields  $ u_{ |_{ \Omega_j } } \ne 0 $ for all $j \in \Upsilon$ and $I'_\Upsilon(u)=0 $. Consequently, by (\ref{tau0}),
$$
\|u\|_{j} > \frac{\tau}{8R} \,\,\, \forall j \in \Upsilon.
$$
This way, $I_\Upsilon (u) \geq c_{\Upsilon}$. But since $ \phi_{ \lambda_n } (u_n) \leq c_\Upsilon $ and $\phi_{ \lambda_n } (u_n) \to I_\Upsilon(u)$,  for $ n $ large, it holds
$$
\|u_n\|_j > \frac{\tau}{8R} \,\,\,  \text{ and } \, \left| \phi_{ \lambda_n}(u_n)-c_\Upsilon \right| \leq \mu, \, \forall j \in \Upsilon.
$$
So $ u_n \in {\cal A}_\mu^{\lambda_n} $, obtaining a contradiction. Thus, we have completed the proof. \qed

In the sequel, $\mu_1,\mu^*$ denote the following numbers
$$
\min_{{\bf t } \in \partial [1/R^2, 1]^l}|I_{\Upsilon}(\gamma_0 ({\bf t }))-c_\Upsilon|=\mu_{1}>0
$$
and
$$
\mu^*=\min\{\mu_1, \delta, {r}/{2}\},
$$
where $\delta$ were given (\ref{A}) and
$$
  r = R^{ 2} \left( \frac{1}{2}-\frac{1}{\theta} \right)^{-1} c_\Upsilon.
$$
Moreover, for each $s>0$, ${ B}_{s}^\lambda$ denotes the set
$$
{ B}_{s}^\lambda = \big\{ u \in E_\lambda \, ; \, \|u\|_{\lambda} \leq s \big\} \,\,\, \text{for} \,\,\, s>0.
$$

\begin{proposition} \label{P}
Let $\mu > 0 $ small enough and $ \Lambda_\ast \ge 1 $ given in the previous proposition. Then, for $ \lambda \ge \Lambda_\ast $, there exists a solution $ u_\lambda $ of $ (A_\lambda) $ such that $ u_\lambda \in {\cal A}_\mu^\lambda \cap \phi_\lambda^{ c_\Upsilon } \cap { B}_{r+1}^\lambda$.

\end{proposition}

\noindent {\bf Proof}
  Let $ \lambda \ge \Lambda_\ast $. Assume that there are no critical points of $ \phi_\lambda $ in $ {\cal A}_\mu^\lambda \cap \phi_\lambda^{ c_\Upsilon } \cap { B}_{r+1}^\lambda $. Since $ \phi_\lambda $ verifies the $ (PS) $ condition, there exists a constant $ d_\lambda > 0 $ such that
$$
   \big\| \phi'_\lambda(u) \big\|_{E^{*}_{\lambda}} \ge d_\lambda, \text{ for all } u \in {\cal A}_\mu^\lambda \cap \phi_\lambda^{ c_\Upsilon } \cap { B}_{r+1}^\lambda.
$$
From Proposition \ref{derivative estimate}, we have
$$
   \big\| \phi'_\lambda(u) \big\|_{E^{*}_{\lambda}} \ge \sigma_0, \text{ for all } u \in \left( {\cal A}_{ 2 \mu }^{\lambda} \setminus {\cal A}_\mu^{\lambda} \right) \cap \phi_{\lambda}^{ c_\Upsilon },
$$
where $ \sigma_0 > 0 $ does not depend on $ \lambda $. In what follows,  $ \Psi \colon E_\lambda \to \mathbb R $ is a continuous functional verifying
$$
   \Psi(u) =  1, \text{ for } u \in {\cal A}_{\frac{3}{2} \mu}^\lambda \cap \Theta_\delta \cap B^{\lambda}_{r},
$$
$$
\ \Psi(u) = 0, \text{ for } u \notin {\cal A}_{2 \mu}^\lambda \cap \Theta_{2\delta} \cap B^{\lambda}_{r+1}
$$
and
$$
0 \le \Psi(u) \le 1, \, \forall u \in E_\lambda.
$$
We also consider $ H \colon \phi_\lambda^{ c_\Upsilon } \to E_\lambda $ given by
$$
   H(u) =
\begin{cases}
   - \Psi(u) \big\| Y(u) \big\|^{ -1 } Y(u), \text{ for } u \in {\cal A}_{2 \mu}^\lambda \cap B^{\lambda}_{r+1}, \\
   \phantom{- \Psi(u) \big\| Y(u) \big\|^{ -1 } Y()} 0, \text{ for } u \notin {\cal A}_{2 \mu}^\lambda \cap B^{\lambda}_{r+1}, \\
\end{cases}
$$
where $ Y $ is a pseudo-gradient vector field for $ \Phi_\lambda $ on $ {\cal K} = \left\{ u \in E_\lambda \, ; \, \phi'_\lambda(u) \ne 0 \right\} $. Observe that  $ H $ is well defined, once $ \phi'_\lambda(u) \ne 0 $, for $ u \in {\cal A}_{2 \mu}^\lambda \cap \phi_\lambda^{ c_\Upsilon } $. The inequality
$$
   \big\| H(u) \big\| \le 1, \, \forall \lambda \ge \Lambda_* \text{ and } u \in \phi_\lambda^{ c_\Upsilon },
$$
guarantees that the deformation flow $ \eta \colon [0, \infty) \times \phi_\lambda^{ c_\Upsilon } \to \phi_\lambda^{ c_\Upsilon } $ defined by
$$
   \frac{d \eta}{dt} = H(\eta), \ \eta(0,u) = u \in \phi_\lambda^{ c_\Upsilon }
$$
verifies
\begin{gather}
   \frac{d}{dt} \phi_\lambda \big( \eta(t,u) \big) \le - \frac{1}{2} \Psi \big( \eta(t,u) \big) \big\| \phi'_\lambda \big( \eta(t,u) \big) \big\| \le 0, \label{eta derivative}\\
	 \left\| \frac{d \eta}{dt} \right\|_\lambda = \big\| H(\eta) \big\|_\lambda \le 1
\end{gather}
and
\begin{equation} \label{eta}
   \eta(t,u) = u \text{ for all } t \ge 0 \text{ and } u \in \phi_\lambda^{ c_\Upsilon } \setminus {\cal A}_{2 \mu}^\lambda \cap B^{\lambda}_{r+1}.
\end{equation}
We study now two paths, which are relevant for what follows: \\

$ \noindent \bullet $ The path $ {\bf t} \mapsto \eta \big( t, \gamma_0( {\bf t} ) \big), \text{ where } \textbf{t} = (t_1,\ldots,t_l) \in [1/R^2, 1]^l $.

\vspace{0.5 cm}

Thereby, if $\mu \in (0,\mu^*)$, we have that
$$
   \gamma_0( {\bf t } ) \notin {\cal A}_{2 \mu}^\lambda, \, \forall {\bf t } \in \partial [1/R^2, 1]^l.
$$
Since
$$
   \phi_\lambda \big( \gamma_0( {\bf t } ) \big) < c_\Upsilon, \, \forall {\bf t } \in \partial [1/R^2, 1]^l,
$$
from (\ref{eta}), it follows that
$$
   \eta \big( t, \gamma_0( {\bf t} ) \big) = \gamma_0( {\bf t} ), \, \forall {\bf t} \in \partial [1/R^2, 1]^l.
$$
So, $ \eta \big( t, \gamma_0( {\bf t} ) \big) \in \Gamma_\ast $, for each $ t \ge 0 $.

\vspace{0.5 cm}

$ \noindent \bullet $ The path $ {\bf t} \mapsto \gamma_0( {\bf t} ), \text{ where } \textbf{t} = (t_1,\ldots,t_l) \in [1/R^2, 1]^l $.

\vspace{0.5 cm}

We observe that
$$
   \text{supp} \big( \gamma_0 ( {\bf t} ) \big)\subset \overline{\Omega_\Upsilon}
$$
and
$$
   \phi_\lambda \big( \gamma_0 ( {\bf t} ) \big) \text{ does not depend on }  \lambda \ge 1,
$$
for all  $ {\bf t} \in [1/R^2, 1]^l $. Moreover,
$$
   \phi_\lambda \big( \gamma_0 ( {\bf t} ) \big) \le c_\Upsilon, \, \forall {\bf t} \in [1/R^2, 1]^l
$$
and
$$
   \phi_\lambda \big( \gamma_0 ( {\bf t} ) \big) = c_\Upsilon \text{ if, and only if, } t_j = \frac{1}{R}, \, \forall j \in \Upsilon.
$$
Therefore
$$
   m_0 = \sup \left\{ \phi_\lambda(u) \, ; \, u \in \gamma_0 \big( [1/R^2,1]^l \big) \setminus A_\mu^\lambda \right\}
$$
is independent of $ \lambda $ and $ m_0 < c_\Upsilon $. Now, observing that there exists $ K_\ast > 0 $ such that
$$
   \big| \phi_{ \lambda }(u) - \phi_{ \lambda}(v) \big| \le K_* \| u-v \|_{ \lambda }, \, \forall u,v \in {\cal B}_r^\lambda ,
$$
we derive
\begin{equation} \label{max estimate}
    \max_{ {\bf t } \in [1/R^2,1]^l } \phi_\lambda \Big( \eta \big( T, \gamma_0 ( {\bf t} ) \big) \Big) \le \max \left\{ m_0, c_\Upsilon-\frac{1}{2 K_\ast} \sigma_0 \mu \right\},
\end{equation}
for $ T > 0 $ large.

In fact, writing $ u = \gamma_0( {\bf t} ) $, $ {\bf t } \in [1/R^2,1]^l $, if $ u \notin A_\mu^\lambda $, from (\ref{eta derivative}),
$$
   \phi_\lambda \big( \eta( t, u ) \big) \le \phi_\lambda (u) \le m_0, \, \forall t \ge 0,
$$
and we have nothing more to do. We assume then $ u \in A_\mu^\lambda $ and set
$$
   \widetilde{\eta}(t) = \eta (t,u), \ \widetilde{d_\lambda} = \min \left\{ d_\lambda, \sigma_0 \right\} \text{ and } T = \frac{\sigma_0 \mu}{K_\ast \widetilde{d_\lambda}}.
$$
Now, we will analyze the ensuing cases: \\

\noindent {\bf Case 1:} $ \widetilde{\eta}(t) \in {\cal A}_{\frac{3}{2} \mu}^\lambda \cap \Theta_\delta \cap B^{\lambda}_{r}, \, \forall t \in [0,T] $.

\noindent {\bf Case 2:} $ \widetilde{\eta}(t_0) \notin {\cal A}_{\frac{3}{2} \mu}^\lambda \cap \Theta_\delta \cap B^{\lambda}_{r}, \text{ for some } t_0 \in [0,T] $. \\

\noindent {\bf Analysis of  Case 1}

In this case, we have $ \Psi \big( \widetilde{\eta}(t) \big) = 1 $ and $ \big\| \phi'_\lambda \big( \widetilde{\eta}(t) \big) \big\| \ge \widetilde{d_\lambda} $ for all $ t \in [0,T] $. Hence, from (\ref{eta derivative}),
$$
   \phi_\lambda \big( \widetilde{\eta}(T) \big) = \phi_\lambda (u) + \int_0^T \frac{d}{ds} \phi_\lambda \big( \widetilde{\eta}(s) \big) \, ds \le c_\Upsilon - \frac{1}{2} \int_0^T \widetilde{d_\lambda} \, ds,
$$
that is,
$$
   \phi_\lambda \big( \widetilde{\eta}(T) \big) \le c_\Upsilon - \frac{1}{2} \widetilde{d_\lambda} T = c_\Upsilon - \frac{1}{2 K_\ast} \sigma_0 \mu,
$$
showing (\ref{max estimate}). \\

\noindent {\bf Analysis of Case 2}:  In this case we have the following situations:
\\

\noindent {\bf (a)}: There exists $t_2 \in [0,T]$ such that $\tilde{\eta}(t_2) \notin \Theta_\delta$, and thus, for $t_1=0$ it follows that
$$
\|\tilde{\eta}(t_2)-\tilde{\eta}(t_1)\| \geq \delta > \mu,
$$
because $\tilde{\eta}(t_1)=u \in \Theta$. \\

\vspace{0.5 cm}

\noindent {\bf (b)}: There exists $t_2 \in [0,T]$ such that $\tilde{\eta}(t_2) \notin B^{\lambda}_r$, so that for $t_1=0$, we get
$$
\|\tilde{\eta}(t_2)-\tilde{\eta}(t_1)\| \geq r > \mu,
$$
because $\tilde{\eta}(t_1)=u \in B^{\lambda}_r$. \\

\noindent {\bf (c)}: \, $\tilde{\eta}(t) \in \Theta_\delta \cap B^{\lambda}_r$ for all $t \in [0,T]$, and there are $0 \leq t_1 \leq t_2 \leq T$ such that $\tilde{\eta}(t) \in  {\cal A}_{\frac{3}{2} \mu}^\lambda \setminus {\cal A}_\mu^\lambda$ for all $t \in [t_1,t_2]$ with
$$
|\phi_\lambda(\tilde{\eta}(t_1))-c_\Upsilon|=\mu \,\,\, \mbox{and} \,\,\, |\phi_\lambda(\tilde{\eta}(t_2))-c_\Upsilon|=\frac{3\mu}{2}
$$
From definition of $K_\ast$, we have
$$
   \| w_2-w_1 \| \ge \frac{1}{K_\ast} \big| \phi_{ \lambda} (w_2) - \phi_{ \lambda} (w_1) \big| \ge \frac{1}{2 K_\ast} \mu.
$$
Then, by  mean value theorem, $ t_2-t_1 \ge \frac{1}{2 K_\ast} \mu $ and, this form,
$$
   \phi_\lambda \big( \widetilde{\eta}(T) \big) \le \phi_\lambda(u) - \int_0^T \Psi \big( \widetilde{\eta}(s) \big) \big\| \phi'_\lambda \big( \widetilde{\eta}(s) \big) \big\| \, ds
$$
implying
$$
   \phi_\lambda \big( \widetilde{\eta}(T) \big) \le c_\Upsilon - \int_{t_1}^{t_2} \sigma_0 \, ds = c_\Upsilon - \sigma_0 (t_2-t_1) \le c_\Upsilon - \frac{1}{2 K_\ast} \sigma_0 \mu,
$$
which proves (\ref{max estimate}). Fixing $ \widehat{\eta} (t_1, \ldots, t_l) = \eta \big( T, \gamma_0 (t_1,\ldots,t_l) \big) $, we have that
$\widehat{\eta}(t_1, \ldots, t_l) \in \Theta_{2\delta}$, and so, $\widehat{\eta}(t_1, \ldots, t_l)|_{\Omega'_j} \not= 0$ for all $j \in \Upsilon$. Thus, $ \widehat{\eta} \in \Gamma_\ast $, leading to
$$
   b_{ \lambda, \Gamma } \le \max_{ (t_1,\ldots,t_l) \in [1/R^2, 1] } \phi_\lambda \big( \widehat{\eta} (t_1,\ldots,t_l) \big) \le \max \left\{ m_0, c_\Upsilon - \frac{1}{2 K_\ast} \sigma_0 \mu \right\} < c_\Upsilon,
$$
which contradicts the fact that $ b_{ \lambda, \Upsilon } \to c_\Upsilon $.
\qed

\vspace{.5cm}
\noindent {\bf [Proof of Theorem \ref{main}]}
According Proposition \ref{P}, for $\mu \in (0, \mu^*)$ and $ \Lambda_\ast \ge 1 $, there exists a solution $ u_\lambda $ for $ (A_\lambda) $ such that $ u_\lambda \in {\cal A}_\mu^\lambda \cap \phi_\lambda^{ c_\Upsilon } $, for all $\lambda \geq \Lambda_*$. \\

\noindent {\bf Claim:}
There are $\lambda_0 \geq \Lambda_*$ and $\mu_0>0$ small enough, such that $u_\lambda$ is a solution for $ \big( P_\lambda \big)$ for $\lambda \geq \Lambda_0$ and $\mu \in (0, \mu_0)$.

Indeed, fixed $\mu \in (0, \mu_0)$, assume by contradiction that there are  $ \lambda_n \to \infty $, such that $(u_{\lambda_n})$ is not a solution for $(SP)_{\lambda_n}$. From Proposition \ref{P}, the sequence $ (u_{\lambda_n}) $ verifies:
\begin{enumerate}
   \item[(a)] $ \phi'_{ \lambda_n }(u_{\lambda_n}) = 0, \, \forall n \in \Bbb N $;
	 \item[(b)] $ \|u_n\|^{2}_{\lambda_n,  \Bbb R^3 \setminus \Omega_\Upsilon }(u_{\lambda_n}) \to 0$;
	 \item[(c)] $ \phi_{ \lambda_n } (u_{\lambda_n}) \to d \leq c_\Upsilon. $
\end{enumerate}	
The item (b) ensures we can use Proposition \ref{P:boundedness of the solutions} to deduce $ u_{\lambda_n} $ is a solution for
$ (SP)_{\lambda_n} $, for large values of $ n $, which is a contradiction, showing this way the claim. \\

Now, our goal is to prove the second part of the theorem. To this end, let  $(u_{\lambda_n})$ be a sequence verifying the above limits. A direct computation gives $ \phi_{ \lambda_n }(u_{ \lambda_n } ) \to d$ with $d \leq c_\Upsilon $. This way, using Proposition \ref{(PS) infty condition} combined with item (c), we derive $( u_{ \lambda_n } )$ converges in $ H^{1}(\Bbb R^3) $ to a function $ u \in H^{1}(\mathbb{R}^3) $, which satisfies $ u = 0 $ outside $ \Omega_\Upsilon $ and $ u_{|_{\Omega_j}} \not= 0, \, j \in \Upsilon $, and $u$ is a positive solution for
         $$
				 \begin{cases}
		        - \Delta u + u +\left(\displaystyle \int_{\Omega_\Upsilon}\frac{u^2}{|x-y|}dy\right)  u= f(u), \text{ in } \Omega_\Upsilon, \\
		        u \in H^{1}_0 ( \Omega_\Upsilon ),
		       \end{cases}
\eqno{(P)_{\infty, \Upsilon}}
$$
and so,
$$
I_{\Upsilon}(u) \geq c_{\Upsilon}.
$$
On the other hand, we also know that
$$
\phi_{ \lambda_n }(u_{ \lambda_n } ) \to I_{\Upsilon}(u),
$$
implying that
$$
I_{\Upsilon}(u)=d \,\,\, \mbox{and} \,\,\, d \geq c_{\Upsilon}.
$$
Since $d \leq c_{\Upsilon}$, we deduce that
$$
I_{\Upsilon}(u)=c_{\Upsilon},
$$
showing that $u$ is a least energy solution for $(P)_{\infty, \Upsilon}$. Consequently, $(u,\phi_u)$ is a least energy solution for the problem	
	$$
      	    \begin{cases}
		     - \Delta u + u + \phi u = f(u), \text{ in } \Omega_\Upsilon, \\
				-\Delta \phi = u^{2}, \,\,\, \text{in} \,\,\mathbb{R}^{3} \\
		     u \in H^{1}_0 ( \Omega_\Upsilon ).
			\end{cases}
  	$$
\qed

\end{document}